# IMPROVING SAMC USING SMOOTHING METHODS: THEORY AND APPLICATIONS TO BAYESIAN MODEL SELECTION PROBLEMS[1]

BY FAMING LIANG

*Texas A&M University*

Stochastic approximation Monte Carlo (SAMC) has recently been proposed by Liang, Liu and Carroll [*J. Amer. Statist. Assoc.* **102** (2007) 305–320] as a general simulation and optimization algorithm. In this paper, we propose to improve its convergence using smoothing methods and discuss the application of the new algorithm to Bayesian model selection problems. The new algorithm is tested through a change-point identification example. The numerical results indicate that the new algorithm can outperform SAMC and reversible jump MCMC significantly for the model selection problems. The new algorithm represents a general form of the stochastic approximation Markov chain Monte Carlo algorithm. It allows multiple samples to be generated at each iteration, and a bias term to be included in the parameter updating step. A rigorous proof for the convergence of the general algorithm is established under verifiable conditions. This paper also provides a framework on how to improve efficiency of Monte Carlo simulations by incorporating some nonparametric techniques.

**1. Introduction.** As known by many researchers, the Metropolis–Hastings (MH) algorithm [Metropolis et al. (1953), Hastings (1970)] and the Gibbs sampler [Geman and Geman (1984)] are prone to get trapped into local energy minima in simulations from a system for which the energy landscape is rugged. In terms of physics, the negative of the logarithmic density/mass function is called the energy function of the system. To overcome the local-trap problem, many advanced Monte Carlo algorithms have been proposed, such as parallel tempering [Geyer (1991), Hukushima and Nemoto (1996)], simulated tempering [Marinari and Parisi (1992), Geyer and Thompson (1995)],

Received September 2006; revised November 2007.
[1]Supported in part by National Science Foundation Grants DMS-04-05748 and DMS-06-07755 and National Cancer Institute Grant CA104620.
*AMS 2000 subject classifications.* 60J22, 65C05.
*Key words and phrases.* Model selection, Markov chain Monte Carlo, reversible jump, smoothing, stochastic approximation Monte Carlo.

evolutionary Monte Carlo [Liang and Wong (2001)], dynamic weighting [Wong and Liang (1997)], multicanonical sampling [Berg and Neuhaus (1991)], $1/k$-ensemble sampling [Hesselbo and Stinchcomble (1995)], the Wang–Landau algorithm [Wang and Landau (2001), Liang (2005)], equi-energy sampler [Mitsutake, Sugita and Okamoto (2003), Kou, Zhou and Wong (2006)], stochastic approximation Monte Carlo (SAMC) [Liang, Liu and Carroll (2007), Atchadé and Liu (2007)], among others. Henceforth, the work by Liang, Liu and Carroll (2007) will be referred to as LLC.

Among the above algorithms, SAMC is a very sophisticated one in both theory and applications. The basic idea of SAMC stems from the Wang–Landau algorithm and can be explained briefly as follows. Let

$$f(x) = c\psi(x), \qquad x \in \mathcal{X}, \tag{1}$$

denote the target probability density/mass function we are working with, where $\mathcal{X}$ is the sample space and $c$ is an unknown constant. Let $E_1, \dots, E_m$ denote a partition of $\mathcal{X}$, and let $\omega_i = \int_{E_i} \psi(x)\,dx$ for $i = 1, \dots, m$. SAMC seeks to sample from the trial distribution

$$f_\omega(x) \propto \sum_{i=1}^{m} \frac{\pi_i \psi(x)}{\omega_i} I(x \in E_i), \tag{2}$$

where $\pi_i$'s are prespecified constants such that $\pi_i > 0$ for all $i$ and $\sum_{i=1} \pi_i = 1$. It is easy to see that if $\omega_1, \dots, \omega_m$ are known, sampling from $f_\omega(x)$ will result in a "random walk" in the space of subregions (by regarding each subregion as a "point") with each subregion being sampled with a frequency proportional to $\pi_i$. Hence, the local-trap problem can be overcome essentially, provided that the sample space is partitioned appropriately. How to partition the sample space will be discussed later. SAMC has been applied successfully to many hard computational problems, such as phylogenetic tree reconstruction [Cheon and Liang (2007)] neural network training [Liang (2007)] and Bayesian model selection [LLC (2007)].

The success of SAMC depends crucially on the estimation of $\omega_i$. LLC propose to estimate $\omega_i$ simultaneously using a stochastic approximation Markov chain Monte Carlo algorithm. Let $\theta_{ti}$ denote the working estimate of $\log(\omega_i/\pi_i)$ obtained at iteration $t$, $\theta_t = (\theta_{t1}, \dots, \theta_{tm})$, and let $\{\gamma_t\}$ be a positive, nonincreasing sequence satisfying the conditions

$$\text{(i)} \quad \sum_{t=1}^{\infty} \gamma_t = \infty, \qquad \text{(ii)} \quad \sum_{t=1}^{\infty} \gamma_t^{\zeta} < \infty, \tag{3}$$

for any $\zeta > 1$. Since $f_\omega(x)$ is invariant to a scale change of $\omega = (\omega_1, \dots, \omega_m)$, that is, $f_{c\omega}(x) = f_\omega(x)$ for any number $c > 0$, the domain of $\theta_t$ can be restricted to a compact set $\Theta$ by adjusting $\theta_t$ with a constant vector, provided that $\Theta$ is large enough. Refer to Chen (2002) and Andrieu, Moulines and Priouret (2005) for more discussions on this issue. The SAMC algorithm iterates between the following two steps.

SAMC ALGORITHM.

(a) Simulate a sample $x_t$ by a single MH update with the invariant distribution

$$f_{\theta_t}(x) \propto \sum_{i=1}^{m} \frac{\psi(x)}{e^{\theta_{ti}}} I(x \in E_i). \tag{4}$$

(b) Set $\theta^* = \theta_t + \gamma_{t+1}(\tilde{\mathbf{e}}_t - \boldsymbol{\pi})$, where $\tilde{\mathbf{e}}_t = (\tilde{e}_{t,1}, \ldots, \tilde{e}_{t,m})$ and $\tilde{e}_{t,i} = 1$ if $x_t \in E_i$ and 0 otherwise. If $\theta^* \in \Theta$, set $\theta_{t+1} = \theta^*$; otherwise, set $\theta_{t+1} = \theta^* + \mathbf{c}^*$, where $\mathbf{c}^* = (c^*, \ldots, c^*)$ can be an arbitrary vector which satisfies the condition $\theta^* + \mathbf{c}^* \in \Theta$.

Under mild conditions, LLC established the convergence of $\theta_t$. The superiority of SAMC in sample space exploration is due to its self-adjusting mechanism. If a subregion is visited, $\theta_t$ will be updated accordingly such that this subregion has a smaller probability to be revisited in the next iteration. Mathematically, if $x_t \in E_i$, then $\theta_{t+1,i} \leftarrow \theta_{t,i} + \gamma_{t+1}(1 - \pi_i)$ and $\theta_{t+1,j} \leftarrow \theta_{t,i} - \gamma_{t+1}\pi_j$ for $j \neq i$. However, this mechanism has not yet reached its maximum efficiency because it does not differentiate between the neighboring and nonneighboring subregions of $E_i$. We note that for many problems, $E_1, \ldots, E_m$ form a sequence of naturally ordered categories with $\omega_1, \ldots, \omega_m$ changing smoothly along the index of subregions. For example, for model selection problems $\mathcal{X}$ can be partitioned according to the index of models, the subregions can be naturally ordered according to the number of parameters contained in each model, and the neighboring subregions often contain similar probability values. Intuitively, $x_t$ may contain some information on its neighboring subregions, so the visiting to its neighboring subregions should also be penalized to some extent in the next iteration. Consequently, this improves the ergodicity of the simulation. Henceforth, we will call a partition with $\omega_1, \ldots, \omega_m$ changing smoothly a smooth partition or say the sample space is partitioned smoothly, and assume that there exists a smooth partition for the problem under study.

In this paper, we show that the efficiency of SAMC can be improved by including at each iteration a smoothing step, which distributes the information contained in each sample to its neighboring subregions. The new algorithm is thus called smoothing-SAMC or SSAMC for simplicity. SSAMC is tested through a change-point identification example in this paper. Our numerical results show that it outperforms both SAMC and reversible jump Markov chain Monte Carlo (RJMCMC) [Green (1995)] for that example. By comparing the sampling mechanisms of SSAMC and RJMCMC, we argue that SSAMC can be superior to RJMCMC for the model selection problems for which the sample space can be partitioned smoothly. A rigorous proof for the convergence of SSAMC is provided in the Appendix. As discussed later,

SSAMC represents the most general form of the stochastic approximation MCMC algorithm.

The remainder of this paper is organized as follows. In Section 2, we describe the SSAMC algorithm and prove a theorem concerning its convergence. In Section 3, we illustrate the use of SSAMC through a mixture Gaussian example. In Section 4, we apply SSAMC to a change-point identification example. In Section 5, we conclude this paper with a brief discussion.

**2. Smoothing-SAMC algorithm.** Suppose that we are working with a distribution as specified in (1) and that the sample space $\mathcal{X}$ has been partitioned into $m$ disjoint subregions $E_1, \ldots, E_m$ according to a function denoted by $\lambda(x)$. Furthermore, we suppose that the subregions have been ordered such that the weights $\omega_1, \ldots, \omega_m$ change smoothly along the index of the subregions.

The SSAMC algorithm is different from the SAMC algorithm in two aspects. First, the gain factor sequence used in SSAMC is a little more restrictive than that used in SAMC. In SSAMC, the gain factor sequence is required to be positive and nonincreasing, and satisfy the following conditions:

$$\text{(5)} \qquad \begin{array}{llll} \text{(i)} & \lim_{t\to\infty} \gamma_t = 0, & \text{(ii)} & \overline{\lim_{t\to\infty}} |\gamma_t^{-1} - \gamma_{t+1}^{-1}| < \infty, \\ \text{(iii)} & \sum_{t=1}^{\infty} \gamma_t = \infty, & \text{(iv)} & \sum_{t=1}^{\infty} \gamma_t^{\zeta} < \infty \qquad \text{for any } \zeta > 1. \end{array}$$

The trade-off is that a higher-order noise term can be included in updating $\theta_t$ as prescribed in (21). In this paper, we set

$$\gamma_t = \frac{T_0}{\max\{T_0, t\}}, \qquad t = 1, 2, \ldots, \tag{6}$$

in all computations, where $T_0$ is a prespecified number. It is easy to see that (6) satisfies the condition (5).

Second, SSAMC allows multiple samples to be generated at each iteration, and employs a smoothed estimate of $p_{ti}$ in updating $\theta_t$, where $p_{ti} = \int_{E_i} f_{\theta_t}(x)\,dx$ is the probability that a sample is drawn from $E_i$ at iteration $t$, and $f_{\theta_t}(x)$ is as defined in (4). Let $x_t^{(1)}, \ldots, x_t^{(\kappa)}$ denote the samples generated by a MH kernel with the invariant distribution $f_{\theta_t}(x)$. Since $\kappa$ is usually a small number, say, 10 to 20, the samples form a sparse frequency vector $\mathbf{e}_{\mathbf{x}_t} = (e_{t1}, \ldots, e_{tm})$ with $e_{ti} = \sum_{l=1}^{\kappa} I(x_t^{(l)} \in E_i)$. Because the law of large numbers does not apply here, $\mathbf{e}_{\mathbf{x}_t}/\kappa$ is not a good estimator of $\mathbf{p}_t = (p_{ti}, \ldots, p_{tm})$. As suggested by many authors, for example, Burman (1987), Hall and Titterington (1987), Dong and Simonoff (1994), Fan, Heckman and Wand (1995) and Aerts, Augustyns and Janssen (1997), the frequency estimate can be improved by a smoothing method. Since we have assumed that the partition is smooth, information in nearby subregions can be borrowed to help produce more accurate estimates of $\mathbf{p}_t$.

In this paper, the frequency estimator $\mathbf{e}_{\mathbf{x}_t}/\kappa$ is smoothed by the Nadaraya–Watson kernel estimator; that is, $p_{ti}$ is estimated by

$$\widehat{p}_{ti} = \frac{\sum_{j=1}^{m} W(\Lambda(i-j)/(mh_t))e_{tj}/\kappa}{\sum_{j=1}^{m} W(\Lambda(i-j)/(mh_t))}, \tag{7}$$

where $W(z)$ is a kernel function with bandwidth $h_t$, and $\Lambda$ is a rough estimate of the range of $\lambda(x)$, $x \in \mathcal{X}$. Here, it is assumed that $W(z)$ has a bounded support; that is, there exists a constant $C$ such that $W(z) = 0$ if $|z| > C$. Under this assumption, it is easy to show that the deviation of $\widehat{p}_{ti}$ from the frequency estimate $e_{ti}/\kappa$ is of the order $O(h_t)$; that is, $\widehat{p}_{ti} - e_{ti}/\kappa = O(h_t)$. Refer to the Appendix [around (32)] for the details of the proof. We have many choices for $W(z)$, for example, an Epanechnikov kernel or a double-truncated Gaussian kernel. The former is standard, and the latter can be written as

$$W(z) = \begin{cases} \exp(-z^2/2), & \text{if } |z| < C, \\ 0, & \text{otherwise.} \end{cases} \tag{8}$$

The bandwidth $h_t$ is chosen as a power function of $\gamma_t$, that is, $h_t = a\gamma_t^b$ for $a > 0$ and $b > 0$. Here $b$ specifies the decay rate of the smoothing adaptation in the SSAMC algorithm. For a small value of $b$, the adaptation can decay very slowly. In all computations of this paper, $W(z)$ is set to the double-truncated Gaussian kernel with $C = 3$ and

$$h_t = \min\left\{\sqrt{\gamma_t}, \frac{\text{range}\{\lambda(x_t^{(1)}), \ldots, \lambda(x_t^{(\kappa)})\}}{2(1 + \log_2(\kappa))}\right\}, \tag{9}$$

where the second term in $\min\{\cdot,\cdot\}$ is the default bandwidth used in conventional density estimation procedures for continuous observations, for example, S-PLUS 5.0 [Venables and Ripley (1999), page 135]. It is easy to see that $h_t = \sqrt{\gamma_t}$ when $t$ becomes large.

In summary, one iteration of the SSAMC algorithm consists of the following three steps:

SSAMC ALGORITHM.

(a) (*Sampling*) Simulate samples $x_t^{(1)}, \ldots, x_t^{(\kappa)}$ using the MH algorithm with the proposal distribution $q(x_t^{(i)}, \cdot)$ and the invariant distribution $f_{\theta_t}(x)$ as defined in (4), where $x_t^{(0)} = x_{t-1}^{(\kappa)}$.

(b) (*Smoothing*) Calculate $\widehat{\mathbf{p}}_t = (\widehat{p}_{t1}, \ldots, \widehat{p}_{tm})$ in (7).

(c) (*Weight updating*) Set

$$\theta^* = \theta_t + \gamma_{t+1}(\widehat{\mathbf{p}}_t - \boldsymbol{\pi}). \tag{10}$$

If $\theta^* \in \Theta$, set $\theta_{t+1} = \theta^*$; otherwise, set $\theta_{t+1} = \theta^* + \mathbf{c}^*$, where $\mathbf{c}^* = (c^*, \dots, c^*)$ can be any vector which satisfies the condition $\theta^* + \mathbf{c}^* \in \Theta$.

For reasons of mathematical convenience, we assume that $\mathcal{X}$ is either finite (for a discrete system) or compact (for a continuum system). For the latter case, $\mathcal{X}$ can be restricted to the region $\{x : \psi(x) \geq \psi_{\min}\}$, where $\psi_{\min}$ is sufficiently small such that the region $\{x : \psi(x) < \psi_{\min}\}$ is not of interest. As in SAMC, $\Theta$ can also be restricted to a compact set. In this paper, we set $\Theta = [-10^{100}, 10^{100}]^m$, although as a practical matter this is essentially equivalent to setting $\Theta = \mathbb{R}^m$. Since both $\mathcal{X}$ and $\Theta$ are compact, it is natural to assume that $f_{\theta_t}(x)$ is bounded away from 0 and $\infty$ on $\mathcal{X}$. Furthermore, we assume that the proposal distribution $q(\cdot,\cdot)$ used in the sampling step of SSAMC satisfies the local positive condition; that is, for every $x \in \mathcal{X}$, there exists $\epsilon_1$ and $\epsilon_2$ such that $q(x,y) \geq \epsilon_2$ if $\|x - y\| \leq \epsilon_1$, where $\|z\|$ denotes the norm of the vector $z$. In a study of MCMC theory, the proposal distribution is often assumed to satisfy the local positive condition [Roberts and Tweedie (1996)].

Under the above assumptions, we establish the following result concerning the convergence of SSAMC. (A formal statement of this result and the proof are given in the Appendix.) As $t \to \infty$, we have

$$\theta_{ti} \to \begin{cases} \text{Const} + \log\left(\int_{E_i} \psi(x)\,dx\right) - \log(\pi_i + \nu), & \text{if } E_i \neq \varnothing, \\ -\infty, & \text{if } E_i = \varnothing, \end{cases} \tag{11}$$

where $\nu = \sum_{j \in \{i : E_i = \varnothing\}} \pi_j/(m - m_0)$ and $m_0$ is the number of empty subregions, and Const represents an arbitrary constant. Since $f_{\theta_t}(x)$ is invariant with respect to a location transformation of $\theta_t$, Const cannot be determined by the samples drawn from $f_{\theta_t}(x)$. To determine the constant term, extra information, for example, $\sum_{i=1}^m e^{\theta_{ti}}$ is equal to a known number, is needed.

LLC discussed several practical issues on implementation of SAMC, including sample space partitioning, convergence diagnostic, and parameter setting (for $\boldsymbol{\pi}$, $T_0$ and the total number of iterations), most of which are still applicable to SSAMC. To make the paper self-contained, they are briefly discussed as follows.

The sample space should be partitioned such that the MH chain can mix reasonably fast within the same subregion. For example, if one chooses to partition the sample space according to the energy function $-\log\psi(x)$, the partition may be done as follows: $E_1 = \{x : -\log\psi(x) \leq u_1\}$, $E_2 = \{x : u_1 < -\log\psi(x) \leq u_2\}, \dots, E_m = \{x : -\log\psi(x) \geq u_{m-1}\}$, with the energy bandwidth $u_i - u_{i-1}$ ($i = 2, \dots, m$) being less than 2, and $u_1$ and $u_m$ being chosen appropriately such that the probabilities contained in $E_1$ and $E_m$ are ignorable. This partition ensures that the MH moves within the same subregion have a reasonable acceptance rate. For the model selection problem, the sample space is usually partitioned according to the model index by assuming that the MH chain can mix reasonably fast in the sample space of each model. If this is not true, one may partition the sample space jointly according to the energy function and the model index.

The convergence of SSAMC can be diagnosed by examining the patterns of the estimates of $\omega$ obtained in multiple runs. If the estimates follow the same pattern, we may reasonably think the runs have been converged. Otherwise, we may think the gain factor is still large at the end of the runs, or some parts of the sample space have not yet been visited, and some parameters should be reset as described below. LLC also proposed to diagnose the convergence of SAMC based on the realized sampling frequencies. This may not work well for SSAMC due to its use of the smoothing estimator at each iteration.

The choice of $\boldsymbol{\pi}$ is problem dependent. If one aims at optimization, $\boldsymbol{\pi}$ may be set biased to low energy regions to improve the ergodicity of the simulation; whereas, if one aims at estimating $\omega$, $\boldsymbol{\pi}$ may be set to a discrete uniform distribution over the subregions. The parameter $T_0$ and the total number of iterations can be determined by a trial and error process based on diagnostics for the convergence of the simulations. If a run is diagnosed as unconverged, SAMC should be rerun with a larger value of $T_0$, a larger number of iterations, or both. In general, a complex problem should associate with a large value of $T_0$ and a large number of iterations.

Below we discuss two more issues specifically related to SSAMC.

- On the choice of smoothing estimators. Theoretically, any smoothing estimator, which satisfies the condition $\widehat{p}_{ti} - e_{ti}/\kappa = O(h^\tau)$ for some $\tau > 0$, can be used in SSAMC. Other than the Nadaraya–Watson kernel estimator, the estimators that possibly can be used include the local log-likelihood estimator [Tibshirani and Hastie (1987), Fan, Heckman and Wand (1995)] and the local polynomial estimator [Aerts, Augustyns and Janssen (1997)], etc. Refer to Simonoff (1998) for a comprehensive review of smoothing estimators.
- On the choice of $\kappa$. Since the convergence of SSAMC is determined by the three parameters $\kappa$, $T_0$ and $N$ (the total number of iterations) together, we suggest that the value of $\kappa$ should be determined together with the values of $T_0$ and $N$ through a trial and error process as described above. In practice, $\kappa$ is usually set to a number less than 20. Since the gain factor is kept at a constant in each iteration, a run with a large $\kappa$ has to end at a large value of $\gamma_t$, provided that the total running time is fixed. The estimates produced by a run ending at a large value of gain factor are often highly variable. In our experience, SSAMC can benefit from the smoothing operation even when $\kappa$ is as small as 5, and the maximum benefit is usually attained at a value of $\kappa$ between 10 and 20.

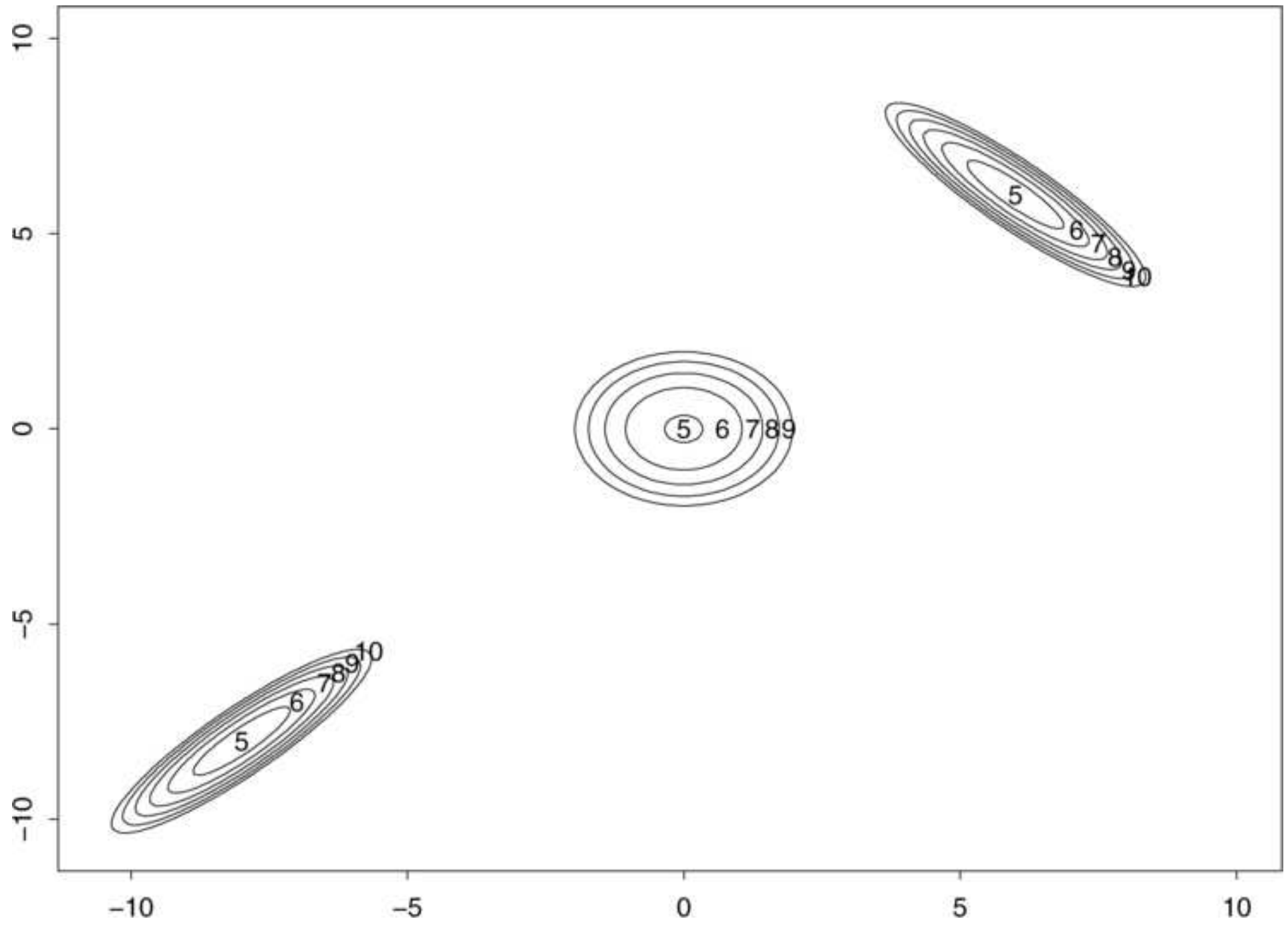

FIG. 1. *Contour plot of the distribution. The numbers in the plot indicate the subregions. For example, $E_5$ includes three separated small areas labeled by "5."*

**3. An illustrative example.** In this section, we illustrate the use of SSAMC through a mixture Gaussian example. Our numerical results indicate that SSAMC can converge much faster than SAMC for this example. Consider the following distribution:

$$f(\mathbf{x}) = \frac{1}{3}N\left[\begin{pmatrix}-8\\-8\end{pmatrix}, \begin{pmatrix}1 & 0.9\\0.9 & 1\end{pmatrix}\right] + \frac{1}{3}N\left[\begin{pmatrix}6\\6\end{pmatrix}, \begin{pmatrix}1 & -0.9\\-0.9 & 1\end{pmatrix}\right] + \frac{1}{3}N\left[\begin{pmatrix}0\\0\end{pmatrix}, \begin{pmatrix}1 & 0\\0 & 1\end{pmatrix}\right],$$

which is identical to an example given in Gilks, Roberts and Sahu (1998), except that the mean vectors are separated by a larger distance in each dimension. Figure 1 shows the contour plot of the distribution, which indicates that the distribution contains three well-separated components. The MH algorithm was first applied to simulate from $f(x)$ with a random walk proposal $N(x, I_2)$, but it failed to mix the three components.

Let $\mathcal{X} = [-10^{100}, 10^{100}]^2$ be compact, and let it be partitioned according to the function $\lambda(x) = -\log f(x)$ into the following subregions: $E_1 = \{x : \lambda(x) < 0.5\}$, $E_2 = \{x : 0.5 \le \lambda(x) < 1.0\}, \ldots, E_{44} = \{x : 21.5 \le \lambda(x) < 22.0\}$ and $E_{45} = \{x : \lambda(x) \ge 22.0\}$, as illustrated by Figure 1. Note that in simulations, we never need to know where the subregions are. It is enough to know which subregion it belongs to for any given sample. SSAMC (also SAMC) offers to learn the weights $\int_{E_1} \psi(x)\,dx/\pi_1, \ldots, \int_{E_m} \psi(x)\,dx/\pi_m$ simultaneously via a stochastic approximation process. The self-adjusting mechanism of SSAMC ensures the success of the learning process; the entire sample space is fully explored and the weights converge to their true values. After convergence, importance samples can then be simulated from the target distribution $f(x)$. Since our purpose of studying this example is just to illustrate how the convergence of SAMC can be accelerated by the smoothing operator, the issue of post-convergence inference will not be discussed here. Refer to LLC for this issue.

SSAMC was run for this example 20 times independently with the setting: $\psi(x) = f(x)$, $T_0 = 25$, $\kappa = 20$, $\Lambda = 22$, $m = 45$, $N = 5 \times 10^5$, and $\pi_1 = \cdots = \pi_m = 1/m$. Table 1 summarizes the estimates of the probabilities $P(E_i) = \int_{E_i} f(x)\,dx$, $i = 5, \ldots, 10$. Note that the subregions $E_1, \ldots, E_4$ are empty. Like SAMC, SSAMC allows the existence of empty subregions in simulations. The corresponding true probability values, which are calculated with a total of $3 \times 10^8$ samples drawn equally from each of the three components of $f(x)$, are also given in Table 1.

For comparison, SAMC was applied to this example with the same setting as that used by SSAMC except for $T_0 = 500$ and $N = 10^7$. SAMC was also run 20 times independently. The computational results are also summarized in Table 1. SSAMC has made a significant improvement over SAMC in terms of accuracy of the estimates; in other words, SSAMC converges much faster than SAMC. On average, the RMSE (root mean squared error) of the SSAMC estimates is only about half of that of the SAMC estimates. Note that under the above settings, we have almost the same gain factor sequence $\{\gamma_t\}$ and exactly the same number of energy evaluations in each of the SAMC and SSAMC runs. Hence, the comparison is fair. The number of energy evaluations is actually a much better measure than the CPU time for comparing efficiency of two algorithms, because the CPU cost is usually dominated by the part used for energy evaluations when we simulate from a complex system, for example, protein folding. We note that this measure has long been used in statistical physics [see, e.g., Hesselbo and Stinchcomble (1995)]. To provide more evidence for the fairness of our comparisons, we also reported in Table 1 the CPU times cost by the above runs; SAMC and SSAMC cost about the same CPU time in each run.

To examine the effect of the sample size $\kappa$ on efficiency of SSAMC, SSAMC was rerun with $\kappa = 5$ and 10. The values of $T_0$ and $N$ were set accordingly,

Table 1
*Comparison of RMSEs (root mean squared errors) of the SAMC and SSAMC estimates*

| | | | SSAMC | | |
|---|---|---|---|---|---|
| **Estimates** | **True prob. (%)** | **SAMC** | $\kappa=5$ | $\kappa=10$ | $\kappa=20$ |
| | 21.60 | 21.65 | 21.66 | 21.68 | |
| $P(E_5)$ | 21.70 | (0.23) | (0.13) | (0.09) | (0.11) |
| | 19.67 | 19.67 | 19.70 | 19.74 | |
| $P(E_6)$ | 19.74 | (0.17) | (0.10) | (0.08) | (0.05) |
| | 23.10 | 23.10 | 23.08 | 23.05 | |
| $P(E_7)$ | 23.04 | (0.18) | (0.09) | (0.08) | (0.07) |
| | 14.01 | 14.01 | 13.99 | 13.98 | |
| $P(E_8)$ | 13.98 | (0.08) | (0.06) | (0.04) | (0.04) |
| | 8.51 | 8.49 | 8.49 | 8.48 | |
| $P(E_9)$ | 8.47 | (0.08) | (0.04) | (0.02) | (0.03) |
| | 5.16 | 5.15 | 5.15 | 5.14 | |
| $P(E_{10})$ | 5.15 | 0.04 | (0.02) | (0.02) | (0.02) |
| CPU (s) | — | 33.2 | 35.6 | 34.8 | 33.9 |

The numbers in the parentheses are the RMSEs of the estimates. The CPU time (in seconds) was measured on a 2.8 GHz computer for a single run of the corresponding algorithm.

$T_0=500/\kappa$ and $N=10^7/\kappa$, such that in these runs we have about the same gain factor sequence and exactly the same number of energy evaluations as in the previous runs. The computational results are also summarized in Table 1. They indicate that the smoothing operation can improve the accuracy of the SAMC estimates generally, even when $\kappa$ is as small as 5.

## 4. Bayesian model selection problems.

LLC applied SAMC to Bayesian model selection problems and compared it to RJMCMC. They conclude that SAMC outperforms RJMCMC when the model space is complex, for example, it contains several modes which are well separated from each other, or some tiny probability models, but, of interest to us. However, when the model space is simple, for example, it only contains several models with comparable probabilities, SAMC may not work better than RJMCMC, as in this case the self-adjusting ability of SAMC is no longer crucial for mixing of the models. In this section, we show that for Bayesian model selection problems, SSAMC can make a significant improvement over SAMC and it can also work better than RJMCMC even when the model space is simple. This is illustrated by a change-point identification example.

The change-point identification problem can be stated as follows. Let $Z=(z_1,z_2,\ldots,z_n)$ denote a sequence of independent observations. Assume that the index set $\{1,2,\ldots,n\}$ has been partitioned into blocks and that the sequence follows the same distribution within blocks; that is, there exists a binary vector $\boldsymbol{\vartheta} = (\vartheta_1, \ldots, \vartheta_{n-1})$ with $\vartheta_{c_1} = \cdots = \vartheta_{c_k} = 1$ and 0 elsewhere, such that

$$0 = c_0 < c_1 < \cdots < c_k < c_{k+1} = n$$

and

$$z_i \sim g_r(\cdot), \qquad c_{r-1} < i \le c_r$$

for $r = 1, 2, \ldots, k+1$, where $g_r(\cdot)$ is a density. Our task is to identify the values of $c_1, \ldots, c_k$.

Recently this problem has been treated by several authors using simulation-based methods, such as the Gibbs sampler [Barry and Hartigan (1993)], jump diffusion [Phillips and Smith (1996)], reversible jump MCMC [Green (1995)] and evolutionary Monte Carlo [Liang and Wong (2000)]. In this article, we follow Barry and Hartigan (1993) to consider the case where $g_r(\cdot)$ is a Gaussian density parameterized by $(\mu_r, \sigma_r^2)$. Let $\boldsymbol{\vartheta}^{(k)}$ denote a configuration of $\boldsymbol{\vartheta}$ with $k$ ones, which represents a model of $k$ change-points. Let $\eta^{(k)} = (\boldsymbol{\vartheta}^{(k)}, \mu_1, \sigma_1^2, \ldots, \mu_{k+1}, \sigma_{k+1}^2)$, $\mathcal{X}_k$ denote the space of models with $k$ change-points, $\boldsymbol{\vartheta}^{(k)} \in \mathcal{X}_k$, and $\mathcal{X} = \bigcup_{k=0}^{n} \mathcal{X}_k$. The log-likelihood of $\eta^{(k)}$ is

$$L(Z|\eta^{(k)}) = -\sum_{i=1}^{k+1} \left\{ \frac{c_i - c_{i-1}}{2} \log \sigma_i^2 + \frac{1}{2\sigma_i^2} \sum_{j=c_{i-1}+1}^{c_i} (z_j - \mu_i)^2 \right\}. \tag{12}$$

To conduct a Bayesian analysis, the following priors are specified for $\eta^{(k)}$. The vector $\boldsymbol{\vartheta}^{(k)}$ is subject to the distribution

$$P(\boldsymbol{\vartheta}^{(k)}) = \frac{\lambda^k}{\sum_{j=0}^{n-1} \lambda^j / j!} \frac{(n-1-k)!}{(n-1)!}, \qquad k = 0, 1, \ldots, n-1,$$

which is equivalent to assuming that $\mathcal{X}_k$ is subject to a truncated Poisson distribution with parameter $\lambda$, and each of the $(n-1)!/[k!(n-1-k)!]$ models in $\mathcal{X}_k$ is a priori equal. The component mean $\mu_i$ is subject to an improper prior, and the component variance $\sigma_i^2$ is subject to an inverse-Gamma $IG(\alpha, \beta)$. By assuming that all the priors are independent, we have the log-prior density,

$$\log P(\eta^{(k)}) = a_k - \sum_{i=1}^{k+1} \left[ (\alpha - 1) \log \sigma_i^2 + \frac{\beta}{\sigma_i^2} \right], \tag{13}$$

where $a_k = (k+1)[\alpha \log \beta - \log \Gamma(\alpha)] + \log(n-1-k)! + k \log \lambda$. The $\alpha$, $\beta$ and $\lambda$ are hyperparameters to be chosen by the user. The log-posterior of $\eta^{(k)}$ (up to an additive constant) can be obtained by adding (12) and (13).

Integrating out the parameters $\mu_1, \sigma_1^2, \ldots, \mu_{k+1}, \sigma_{k+1}^2$ from the full posterior distribution and taking a logarithm, we have

$$
\begin{aligned}
\log P(\boldsymbol{\vartheta}^{(k)}|Z) = a_k &+ \frac{k+1}{2}\log 2\pi \\
&- \sum_{i=1}^{k+1}\Biggl\{\frac{1}{2}\log(c_i - c_{i-1}) - \log\Gamma\left(\frac{c_i - c_{i-1} - 1}{2} + \alpha\right) \\
&+ \left(\frac{c_i - c_{i-1} - 1}{2} + \alpha\right) \\
&\times \log\left[\beta + \frac{1}{2}\sum_{j=c_{i-1}+1}^{c_i} z_j^2 - \frac{(\sum_{j=c_{i-1}+1}^{c_i} z_j)^2}{2(c_i - c_{i-1})}\right]\Biggr\}.
\end{aligned}
\tag{14}
$$

The MAP (maximum a posteriori) estimate of $\boldsymbol{\vartheta}$ is often a reasonable solution to the problem. In practice, we are also interested in estimating the marginal posterior distribution $P(\mathcal{X}_k|Z)$. SSAMC can be applied to estimate this distribution. Without loss of generality, we restrict our consideration to the models with $k_{\min} \le k \le k_{\max}$. Let $E_k = \mathcal{X}_k$ and $\psi(\cdot) \propto P(\boldsymbol{\vartheta}^{(k)}|Z)$. It follows from (11) that $\widehat{\omega}_i^{(t)}/\widehat{\omega}_j^{(t)} = e^{\theta_{ti} - \theta_{tj}}$ forms a consistent estimator for the ratio $P(\mathcal{X}_i|Z)/P(\mathcal{X}_j|Z)$ when $t$ is large.

For the change-point identification problem, the sampling step of SSAMC can be performed as follows. Let $\boldsymbol{\vartheta}_t^{(k,l)}$ denote the $l$th sample generated at iteration $t$, where $k$ indicates the number of change-points of the sample. The next sample can be generated according to the following procedure:

(a) Set $j = k-1$, $k$, or $k+1$ according to probabilities $q_{k,j}$, where $q_{k,k} = \frac{1}{3}$ for $k_{\min} \le k \le k_{\max}$, $q_{k_{\min},k_{\min}+1} = q_{k_{\max},k_{\max}-1} = \frac{2}{3}$, and $q_{k,k+1} = q_{k,k-1} = \frac{1}{3}$ if $k_{\min} < k < k_{\max}$.

(b) If $j = k$, update $\boldsymbol{\vartheta}_t^{(k,l)}$ by a "simultaneous" move (described below); if $j = k+1$, update $\boldsymbol{\vartheta}_t^{(k,l)}$ by a "birth" move (described below); and if $j = k-1$, update $\boldsymbol{\vartheta}_t^{(k,l)}$ by a "death" move (described below).

The "birth," "death" and "simultaneous" moves are designed similarly to those described in Green (1995). In the "birth" move, a random number, say $u$, is first drawn uniformly from the set $\{0, 1, \ldots, k\}$; then another random number, say $v$, is drawn uniformly from the set $\{c_u + 1, \ldots, c_{u+1} - 1\}$, and it is proposed to set $\vartheta_v = 1$. The resulting new sample is denoted by $\boldsymbol{\vartheta}_*^{(k+1)}$. In the "death" move, a random number, say $u$, is drawn uniformly from the set $\{1, 2, \ldots, k\}$, and it is proposed to set $\vartheta_{c_u} = 0$. The resulting new sample is denoted by $\boldsymbol{\vartheta}_*^{(k-1)}$. In the "simultaneous" move, a random number, say $u$, is first randomly drawn from the set $\{1, 2, \ldots, k\}$; then another random number, say $v$, is uniformly drawn from the set $\{c_{u-1} + 1, \ldots, c_u - 1, c_u +$ $1,\dots,c_{u+1}-1\}$, and it is proposed to set $\vartheta_{c_u}=0$ and $\vartheta_v=1$. The resulting new sample is denoted by $\boldsymbol{\vartheta}_*^{(k)}$. The acceptance probabilities of the three types of moves are as follows. For the "birth" move, it is

$$\min\left\{1, \frac{e^{\theta_{tk}}}{e^{\theta_{t,k+1}}} \frac{P(\boldsymbol{\vartheta}_*^{(k+1)}|X)}{P(\boldsymbol{\vartheta}_t^{(k,l)}|X)} \frac{q_{k+1,k}}{q_{k,k+1}} \frac{c_{u+1}-c_u-1}{1}\right\}. \tag{15}$$

For the "death" move, it is

$$\min\left\{1, \frac{e^{\theta_{tk}}}{e^{\theta_{t,k-1}}} \frac{P(\boldsymbol{\vartheta}_*^{(k-1)}|X)}{P(\boldsymbol{\vartheta}_t^{(k,l)}|X)} \frac{q_{k-1,k}}{q_{k,k-1}} \frac{1}{c_{u+1}-c_{u-1}-1}\right\}. \tag{16}$$

For the "simultaneous" move, it is

$$\min\left\{1, \frac{P(\boldsymbol{\vartheta}_*^{(k)}|X)}{P(\boldsymbol{\vartheta}_t^{(k,l)}|X)}\right\}, \tag{17}$$

because the proposal densities are symmetric in the sense $T(\boldsymbol{\vartheta}_t^{(k,l)} \to \boldsymbol{\vartheta}_*^{(k)}) = T(\boldsymbol{\vartheta}_*^{(k)} \to \boldsymbol{\vartheta}_t^{(k,l)}) = 1/(c_{u+1}-c_{u-1}-2)$.

Our simulated dataset consists of 1000 observations with $z_1,\dots,z_{120} \sim N(-0.5,1)$, $z_{121},\dots,z_{210} \sim N(0.5,0.5)$, $z_{211},\dots,z_{460} \sim N(0,1.5)$, $z_{461},\dots,z_{530} \sim N(-1,1)$, $z_{531},\dots,z_{615} \sim N(0.5,2)$, $z_{616},\dots,z_{710} \sim N(1,1)$, $z_{711},\dots,z_{800} \sim N(0,1)$, $z_{801},\dots,z_{950} \sim N(0.5,0.5)$ and $z_{951},\dots,z_{1000} \sim N(1,1)$. The time plot is shown in Figure 2. For this dataset we set the hyperparameters $\alpha=\beta=0.05$ and $\lambda=1$. In simulations, we set $k_{\min}=7$ and $k_{\max}=14$. The values of $k_{\min}$ and $k_{\max}$ can be determined rapidly with a short pilot run of the above algorithm. Outside this range, we have $P(\mathcal{X}_i|Z)\approx 0$. SSAMC was run 20 times independently with $\kappa=20$, $T_0=5$, $N=10^5$, $\Lambda=k_{\max}-k_{\min}+1$, $m=8$, and $\pi_1=\cdots=\pi_m=\frac{1}{m}$. The results are summarized in Figure 2 and Table 2.

Figure 2 compares the true change-point pattern and its MAP estimate, which are $(120,210,460,530,615,710,800,950)$ and $(120,211,460,531,610,709,801,939)$, respectively. The largest discrepancy of the two patterns occurs at the last change-point position. A detailed exploration of the original data gives a strong support to the MAP estimate. The last ten observations of the second last cluster have a larger mean value than the expected and thus, they tend to be clustered to the last cluster. Our computation shows that the log-posterior probability of the MAP estimate is 5.33 higher than that of the true pattern.

For comparison, SAMC and RJMCMC were also applied to this example. Each algorithm was run 20 times independently. The computational results are summarized in Table 2. SAMC employs the same sample space partition, the same transition proposals (in the sampling step) and the same

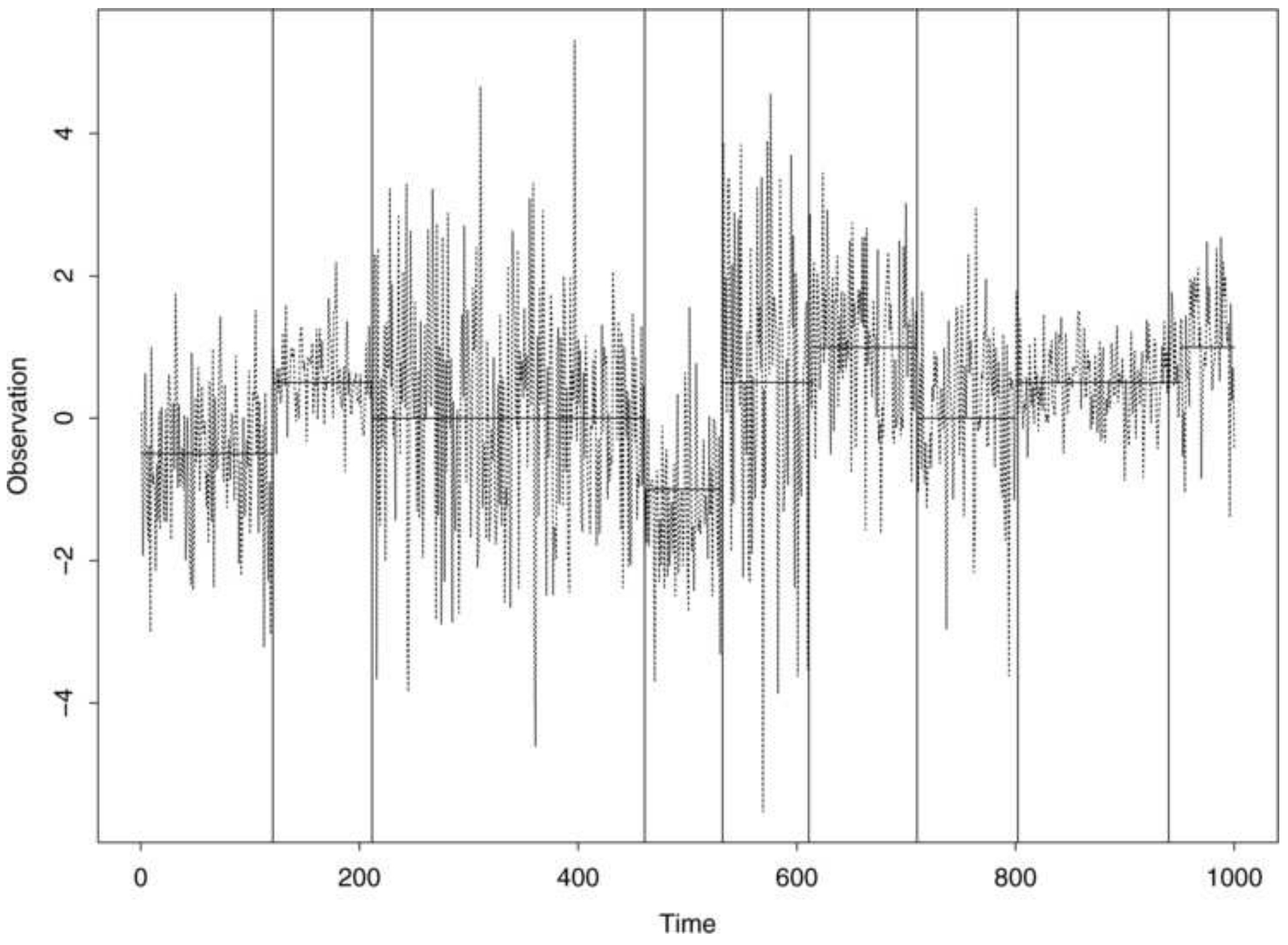

FIG. 2. *Comparison of the true change-point pattern (horizontal lines) and its MAP estimate (vertical lines).*

parameter setting as SSAMC except for $T_0 = 100$ and $N = 2 \times 10^6$. RJMCMC employs the same transition proposals as those used by SSAMC and

TABLE 2
*The estimated posterior distribution $P(\mathcal{X}_k|Z)$ for the change-point identification example*

| | **SSAMC** | | **SAMC** | | **MSAMC** | | **RJMCMC** | |
|---|---|---|---|---|---|---|---|---|
| ***k*** | **prob(%)** | **SD** | **prob(%)** | **SD** | **prob(%)** | **SD** | **prob(%)** | **SD** |
| 7 | 0.1010 | 0.0023 | 0.0944 | 0.0029 | 0.0978 | 0.0020 | 0.0907 | 0.0046 |
| 8 | 55.4666 | 0.2470 | 55.3928 | 0.6112 | 55.0810 | 0.3507 | 55.5726 | 0.3451 |
| 9 | 33.3744 | 0.1659 | 33.3728 | 0.3573 | 33.3798 | 0.2228 | 33.2117 | 0.2052 |
| 10 | 9.2982 | 0.1026 | 9.3647 | 0.2788 | 9.5903 | 0.1351 | 9.3537 | 0.1441 |
| 11 | 1.5655 | 0.0287 | 1.5785 | 0.0685 | 1.6457 | 0.0304 | 1.5694 | 0.0400 |
| 12 | 0.1768 | 0.0042 | 0.1803 | 0.0097 | 0.1871 | 0.0042 | 0.1845 | 0.0097 |
| 13 | 0.0157 | 0.0005 | 0.0154 | 0.0009 | 0.0166 | 0.0004 | 0.0165 | 0.0011 |
| 14 | 0.0018 | 0.0001 | 0.0011 | 0.0001 | 0.0018 | 0.0001 | 0.0009 | 0.0002 |
| CPU (s) | 25.8 | | 25.5 | | 24.9 | | 23.9 | |

SD: standard deviation of the estimates. CPU: the CPU time (in seconds) cost by a single run of the corresponding algorithm on a 2.8 GHz computer.

SAMC, and performs $2 \times 10^6$ iterations in each run. Therefore, in each run the three algorithms perform exactly the same number of energy evaluations, and SAMC and SSAMC also employ the same gain factor sequence. The comparisons made in Table 2 are thus fair to each of the algorithms.

SSAMC works best for this example among the three algorithms. As known by us, RJMCMC can be regarded as a general MH algorithm; and for such a simple problem, it is really hard to find another Monte Carlo algorithm to beat it. However, SSAMC does. SSAMC is different from RJMCMC in two respects. First, like SAMC, it has the capability to self-adjust the acceptance rate of the moves. This capability enables it to overcome any difficulties in the dimension jumping moves and to explore the entire model space very quickly. Second, it has the capability to make use of nearby model information to improve its estimation. However, this can hardly be done in RJMCMC due to the strict requirement for its Markovian property. These two capabilities make SSAMC potentially more efficient than RJMCMC for all types of Bayesian model selection problems. This is evidenced by the observations: LLC showed that SAMC can make a significant improvement over RJMCMC for complex Bayesian model selection problems, and in this paper we showed that SSAMC can make a significant improvement over RJMCMC for simple Bayesian model selection problems.

It is worth pointing out that although the overall performance of SAMC is worse than that of RJMCMC for this example, SAMC tends to work better than RJMCMC for the low probability model spaces, for example, the spaces with 7 and 14 change-points. This is due to the fact that SAMC samples equally from each model space, while RJMCMC samples from each model space proportionally to its probability.

We have also tried a variation of SSAMC for this example. At each step, only multiple samples are generated, but no smoothing is operated on the frequency estimator; that is, (10) in the SSAMC algorithm is replaced by (18),

$$\theta^* = \theta_t + \gamma_{t+1}(\mathbf{e}_{\mathbf{x}_t}/\kappa - \boldsymbol{\pi}). \tag{18}$$

It was run 20 times with exactly the same setting as that used by SSAMC. The results reported in Table 2 under the column MSAMC indicate that averaging over multiple samples can improve the convergence of SAMC, but a further smoothing operation on the frequency estimator is also important.

Later, SSAMC was rerun with some other settings, for example, $\kappa = 10$, $T_0 = 10$ and $N = 2 \times 10^5$; it yielded similar results to those reported in Table 2.

**5. Discussion.** In this paper, we have introduced the SSAMC algorithm, studied its convergence and discussed its application to Bayesian model selection problems. Our numerical results show that SSAMC can converge much faster than SAMC when the sample space is partitioned smoothly. For the problems for which the partition contains abrupt jumps between neighboring subregions, smoothing could potentially make the estimation worse. In the latter case, the subregions can be reordered according to the estimates of $\omega_i$'s from a pilot run such that the resulting partition is smooth.

This paper has made two contributions to the literature. First, it establishes the convergence of a stochastic approximation Markov chain Monte Carlo algorithm under verifiable conditions. The algorithm we studied is very general. It allows multiple samples to be generated at each iteration and a higher-order bias term to be included in the weight updating step. The existing stochastic approximation MCMC algorithms usually only allow a single sample to be generated at each iteration [e.g., Benveniste, Métivier and Priouret (1990), Tadić (1997) and Andrieu, Moulines and Priouret (2005)]. Younes (1999) proved convergence for a stochastic approximation MCMC algorithm which allows for multiple samples, but does not allow for the higher-order bias term. In addition, the conditions assumed by Younes (1999) are less verifiable than those assumed in this paper. Gu and Kong (1998) also studied convergence of a stochastic approximation MCMC algorithm under less verifiable conditions. As indicated by Benveniste, Métivier and Priouret (1990), Chen (2002) and Kushner and Ying (2003), the convergence theory established in this paper can be applied in a much broader context, such as signal processing and adaptive control.

Second, this paper shows how to improve SAMC using smoothing methods, which also provides a general framework on how to improve efficiency of Monte Carlo simulations by incorporating nonparametric techniques. For an illustrative purpose, we employ the Nadaraya–Watson kernel estimator. An advanced smoothing technique, such as the local log-likelihood estimator, should work better in general. It is worth noting that even when the smoothing adaptation, which can decay extremely slowly by choosing $h_t = a\gamma_t^b$ and $b$ a small positive number, stops, SSAMC does not become the same as SAMC. SAMC only allows a single sample to be generated in each iteration, while SSAMC allows multiple samples to be generated in each iteration. Also, the theory established for SAMC by LLC cannot be directly extended to the case of multiple samples. Allowing multiple samples to be generated in each iteration is important, as it provides us much freedom to incorporate some data-mining techniques into simulations. We hope the present work will trigger more research in this direction.

LLC discussed the applications of SAMC to the problems for which the sample space is jointly partitioned according to two functions $\lambda_1(x)$ and $\lambda_2(x)$. The applicability of SSAMC to these problems is apparent. For the joint partitions, the subregions can usually be naturally ordered as a contingency table, and the smoothing estimator used in this paper can be easily extended to the table. Our preliminary results (not reported here) show that in this case, the superiority of SSAMC over SAMC is even more significant.

## APPENDIX: THEORETICAL RESULTS ON SAMC

The Appendix is organized as follows. In Section A.1, we describe a theorem for the convergence of the SSAMC algorithm. In Section A.2, we describe a general version of the SSAMC algorithm and give conditions for its convergence. In Section A.3, we prove the convergence of the general algorithm described in Section A.2. In Section A.4, we prove the convergence of the SSAMC algorithm by verifying that it satisfies the convergence conditions of the general algorithm.

**A.1. A convergence theorem for SSAMC.** Without loss of generality, we only show the convergence presented in (11) for the case that all subregions are nonempty, that is, $\nu = 0$. Extension to the case $\nu \neq 0$ is trivial, since replacing (10) by (19) (given below) will not change the process of SSAMC simulation:

$$\theta' = \theta_t + \gamma_t(\widehat{\mathbf{p}}_{t+1} - \boldsymbol{\pi} - \boldsymbol{\nu}), \tag{19}$$

where $\boldsymbol{\nu} = (\nu, \dots, \nu)$ is an $m$-vector of $\nu$.

THEOREM A.1. *Let $E_1, \dots, E_m$ be a partition of a compact sample space $\mathcal{X}$ and $\psi(x)$ be a nonnegative function defined on $\mathcal{X}$ with $0 < \int_{E_i} \psi(x)\,dx < \infty$ for all $E_i$'s. Let $\boldsymbol{\pi} = (\pi_1, \dots, \pi_m)$ be an $m$-vector with $0 < \pi_i < 1$ and $\sum_{i=1}^m \pi_i = 1$. Let $\Theta$ be a compact set of $m$ dimensions, and there exists a constant $C$ such that $\breve{\theta} \in \Theta$, where $\breve{\theta} = (\breve{\theta}_1, \dots, \breve{\theta}_m)$ and $\breve{\theta}_i = C + \log(\int_{E_i} \psi(x)\,dx) - \log(\pi_i)$. Let $\theta_0 \in \Theta$ be an initial estimate of $\breve{\theta}$, and $\theta_t \in \Theta$ be the estimate of $\breve{\theta}$ at iteration $t$. Let $\{\gamma_t\}$ be a nonincreasing, positive sequence satisfying (5). Let the bandwidth $h_t$ be a power function of $\gamma_t$, that is, $h_t = a\gamma_t^b$ for some $a > 0$ and $b > 0$, when $t$ becomes large. Suppose that $f_{\theta_t}(x)$ is bounded away from 0 and $\infty$ on $\mathcal{X}$, and the proposal distribution satisfies the local positive condition. As $t \to \infty$, we have*

$$P\left\{\lim_{t\to\infty} \theta_{ti} = \mathrm{Const} + \log\left(\int_{E_i} \psi(x)\,dx\right) - \log(\pi_i)\right\} = 1, \qquad i = 1, \dots, m, \tag{20}$$

*where* Const *represents an arbitrary constant.*

**A.2. A convergence theorem for a general stochastic approximation algorithm.** Let $\mathbf{x}_t = (x_t^{(1)}, \ldots, x_t^{(\kappa)})$ be the collection of the samples generated by a MH kernel at iteration $t$, $\mathbf{f}_{\theta_t}(\mathbf{x})$ be the invariant distribution of the MH kernel, $H(\theta_t, \mathbf{x}_{t+1}) = \mathbf{e}_{\mathbf{x}_{t+1}}/\kappa - \boldsymbol{\pi}$, $h(\theta) = \int_{\mathcal{X}^k} H(\theta, \mathbf{x})\mathbf{f}_\theta(d\mathbf{x})$, and $\xi_{t+1} = H(\theta_t, \mathbf{x}_{t+1}) - h(\theta_t) + \gamma_{t+1}^\tau \eta(\mathbf{x}_{t+1})$ with $\tau > 0$ and $\eta(\mathbf{x}_{t+1})$ being a bounded function of $\mathbf{x}_{t+1}$; that is, there exists a constant $\Delta$ such that $\|\eta(\mathbf{x}_{t+1})\| \leq \Delta$ for all $t = 0, 1, \ldots$. The SSAMC algorithm can then be expressed in a more general form by replacing (10) by (21):

$$\theta^* = \theta_t + \gamma_{t+1} h(\theta_t) + \gamma_{t+1}\xi_{t+1}. \tag{21}$$

In the following the general stochastic approximation MCMC algorithm is analyzed under the following conditions.

*Conditions on the step-sizes.*

($A_1$) The sequence $\{\gamma_t\}_{t=0}^\infty$ is nonincreasing, positive and satisfies the condition (5).

*Drift conditions on the transition kernel $P_\theta$.* Below, we first give some definitions on general drift and continuity conditions, and then give the specific drift and continuity conditions for the SSAMC algorithm.

Assume that a transition kernel $P$ is $\psi$-irreducible, aperiodic and has a stationary distribution on a sample space denoted by $\mathcal{X}$. A set $C \subset \mathcal{X}$ is said to be small if there exists a probability measure $\nu$ on $\mathcal{X}$, a positive integer $l$ and $\delta > 0$ such that

$$P_\theta^l(x, A) \geq \delta\nu(A) \qquad \forall x \in C,\ \forall A \in \mathcal{B}_{\mathcal{X}},$$

where $\mathcal{B}_{\mathcal{X}}$ is the Borel set of $\mathcal{X}$. A function $V : \mathcal{X} \to [1, \infty)$ is said to be a drift function outside $C$ if there exist constants $\lambda < 1$ and $b$ such that

$$P_\theta V(x) \leq \lambda V(x) + bI(x \in C) \qquad \forall x \in \mathcal{X},$$

where $P_\theta V(x) = \int_{\mathcal{X}} P_\theta(x, \mathbf{y})V(\mathbf{y})\,d\mathbf{y}$. For $g : \mathcal{X} \to \mathbb{R}^d$, define the norm

$$\|g\|_V = \sup_{x \in \mathcal{X}} \frac{|g(x)|}{V(x)},$$

and define the set $\mathcal{L}_V = \{g : \mathcal{X} \to \mathbb{R}^d, \|g\|_V < \infty\}$.

The specific drift and continuity conditions for the SSAMC algorithm can be described as follows. Let $\mathbf{P}_{\boldsymbol{\theta}}$ be the joint transition kernel for generating the samples $\mathbf{x} = (x^{(1)}, \ldots, x^{(\kappa)})$ at each iteration by ignoring the subscript $t$, $\mathcal{X}^\kappa = \mathcal{X} \times \cdots \times \mathcal{X}$ be the product sample space, $\mathbf{A} = A_1 \times \cdots \times A_\kappa$ be a measurable rectangle in $\mathcal{X}^\kappa$ for which $A_i \in \mathcal{B}_{\mathcal{X}}$ for $i = 1, \ldots, \kappa$, and $\mathcal{B}_{\mathcal{X}^\kappa} = \mathcal{B}_{\mathcal{X}} \times \cdots \times \mathcal{B}_{\mathcal{X}}$ be the $\sigma$-algebra generated by measurable rectangles.

($A_2$) The transition kernel $\mathbf{P}_\theta$ is irreducible and aperiodic for any $\theta \in \Theta$. There exist a function $V : \mathcal{X}^\kappa \to [1, \infty)$ and constants $\alpha \geq 2$ and $\beta \in (0, 1]$ such that:

(i) For any $\theta \in \Theta$, there exist a set $\mathbf{C} \subset \mathcal{X}^\kappa$, an integer $l$, constants $0 < \lambda < 1$, $b, \varsigma$, $\delta > 0$ and a probability measure $\nu$ such that

$$\mathbf{P}_\theta^l V^\alpha(\mathbf{x}) \le \lambda V^\alpha(\mathbf{x}) + bI(\mathbf{x} \in \mathbf{C}) \qquad \forall \mathbf{x} \in \mathcal{X}^\kappa, \tag{22}$$

$$\mathbf{P}_\theta V^\alpha(\mathbf{x}) \le \varsigma V^\alpha(\mathbf{x}) \qquad \forall \mathbf{x} \in \mathcal{X}^\kappa, \tag{23}$$

$$\mathbf{P}_\theta^l(\mathbf{x}, \mathbf{A}) \ge \delta \nu(\mathbf{A}) \qquad \forall \mathbf{x} \in \mathbf{C},\ \forall \mathbf{A} \in \mathcal{B}_{\mathcal{X}^\kappa}. \tag{24}$$

(ii) There exists a constant $c_1$ such that for all $\mathbf{x} \in \mathcal{X}^\kappa$ and $\theta, \theta' \in \Theta$,

$$\|H(\theta, \mathbf{x})\| \le c_1 V(\mathbf{x}), \tag{25}$$

$$\|H(\theta, \mathbf{x}) - H(\theta', \mathbf{x})\| \le c_1 V(\mathbf{x}) \|\theta - \theta'\|^\beta. \tag{26}$$

(iii) There exists a constant $c_2$ such that for all $\theta, \theta' \in \Theta$,

$$\|\mathbf{P}_\theta g - \mathbf{P}_{\theta'} g\|_V \le c_2 \|g\|_V |\theta - \theta'|^\beta \qquad \forall g \in \mathcal{L}_V, \tag{27}$$

$$\|\mathbf{P}_\theta g - \mathbf{P}_{\theta'} g\|_{V^\alpha} \le c_2 \|g\|_{V^\alpha} |\theta - \theta'|^\beta \qquad \forall g \in \mathcal{L}_{V^\alpha}. \tag{28}$$

*Lyapunov condition on $h(\theta)$.* Let $\mathcal{L} = \{\theta \in \Theta : h(\theta) = \mathbf{0}\}$.

($A_3$) The function $h : \Theta \to \mathbb{R}^d$ is continuous, and there exists a continuously differentiable function $v : \Theta \to [0, \infty)$ such that $\dot{v}(\theta) = \nabla^T v(\theta) h(\theta) < 0$, $\forall \theta \in \mathcal{L}^c$ and $\sup_{\theta \in Q} \dot{v}(\theta) < 0$ for any compact set $Q \subset \mathcal{L}^c$.

*A main convergence result.* Let $\mathbb{P}_{\mathbf{x}_0, \theta_0}$ denote the probability measure of the Markov chain $\{(\mathbf{x}_t, \theta_t)\}$, started in $(\mathbf{x}_0, \theta_0)$, and implicitly defined by the sequences $\{\gamma_t\}$. Also define $D(\mathbf{z}, A) = \inf_{\mathbf{z}' \in A} \|\mathbf{z} - \mathbf{z}'\|$.

THEOREM A.2. *Assume the conditions* ($A_1$), ($A_2$) *and* ($A_3$) *hold, and* $\sup_{\mathbf{x} \in \mathcal{X}^\kappa} V(\mathbf{x}) < \infty$. *Let the sequence* $\{\theta_n\}$ *be defined as in the stochastic approximation algorithm. Then for all* $(\mathbf{x}_0, \theta_0) \in \mathcal{X}^\kappa \times \Theta$,

$$\lim_{t \to \infty} D(\theta_t, \mathcal{L}) = 0, \qquad \mathbb{P}_{\mathbf{x}_0, \theta_0}\text{-}a.e.$$

**A.3. Proof of Theorem A.2.** The following lemma is a partial restatement of Proposition 6.1 of Andrieu, Moulines and Priouret (2005).

LEMMA A.1. *Assume the drift condition* ($A_2$). *Then the following results hold:*

($B_1$) *For any $\theta \in \Theta$, the Markov kernel $\mathbf{P}_\theta$ has a single stationary distribution $\mathbf{f}_\theta$. In addition $H : \Theta \times \mathcal{X}^\kappa$ is measurable for all $\theta \in \Theta$, $h(\theta) = \int_{\mathcal{X}^\kappa} H(\theta, \mathbf{x}) \mathbf{f}_\theta(d\mathbf{x}) < \infty$.*

($B_2$) *For any $\theta \in \Theta$, the Poisson equation $u(\theta, \mathbf{x}) - \mathbf{P}_\theta u(\theta, \mathbf{x}) = H(\theta, \mathbf{x}) - h(\theta)$ has a solution $u(\theta, \mathbf{x})$, where $\mathbf{P}_\theta u(\theta, \mathbf{x}) = \int_{\mathcal{X}^\kappa} u(\theta, \mathbf{x}') \mathbf{P}_\theta(\mathbf{x}, \mathbf{x}')\, d\mathbf{x}'$. There exist a function $V : \mathcal{X}^\kappa \to [1, \infty)$ such that the set $\{\mathbf{x} \in \mathcal{X}^\kappa : V(\mathbf{x}) < \infty\} \ne \varnothing$, constant $\beta \in (0, 1]$, $p \ge 2$ such that for any compact subset $\Theta_0 \subset \Theta$,*

(i) $\sup_{\theta\in\Theta_0}\|H(\theta,\mathbf{x})\|_V<\infty$,

(ii) $\sup_{\theta\in\Theta_0}(\|u(\theta,\mathbf{x})\|_V+\|\mathbf{P}_\theta u(\theta,\mathbf{x})\|_V)<\infty$,

(iii) $\sup_{(\theta,\theta')\in\Theta_0}|\theta-\theta'|^{-\beta}(\|u(\theta,\mathbf{x})-u(\theta,\mathbf{x}')\|_V+\|\mathbf{P}_\theta u(\theta,\mathbf{x})-\mathbf{P}_{\theta'}u(\theta', \mathbf{x})\|_V)<\infty$.

Vladislav Tadić studied the convergence of a stochastic approximation MCMC algorithm, which is the same as the SAMC algorithm except that it does not include the step of sample space partitioning, under different conditions from those given in Andrieu, Moulines and Priouret (2005). Tadić proved the following lemma, which corresponds to Theorem 4.1 and Lemma 2.2 of Tadić (1997). In this paper, we show that the results also hold for SSAMC. Our proof is similar to Tadić's except for some necessary changes for including the higher-order noise term in $\xi_t$.

LEMMA A.2. *Assume that the conditions* ($A_1$), ($A_2$), ($B_1$) *and* ($B_2$) *hold and that* $\sup_{\mathbf{x}\in\mathcal{X}^\kappa}V(\mathbf{x})<\infty$. *For the SSAMC algorithm, the following results hold:*

($C_1$) *There exist* $\mathbb{R}^d$*-valued random processes* $\{\epsilon_t\}_{t\ge0}$, $\{\epsilon'_t\}_{t\ge0}$ *and* $\{\epsilon''_t\}_{t\ge0}$ *defined on a probability space* $(\Omega,\mathcal{F},\mathbb{P})$ *such that*

$$\gamma_{t+1}\xi_{t+1}=\epsilon_{t+1}+\epsilon'_{t+1}+\epsilon''_{t+1}-\epsilon''_t,\qquad t\ge0. \tag{29}$$

($C_2$) *The series* $\sum_{t=0}^\infty\|\epsilon'_t\|$, $\sum_{t=0}^\infty\|\epsilon''_t\|^2$ *and* $\sum_{t=0}^\infty\|\epsilon_{t+1}\|^2$ *all converge a.s. and*

$$E(\epsilon_{t+1}|\mathcal{F}_t)=0,\qquad a.s.,\ n\ge0, \tag{30}$$

*where* $\{\mathcal{F}_t\}_{t\ge0}$ *is a family of* $\sigma$*-algebras of* $\mathcal{F}$ *satisfying* $\sigma\{\theta_0\}\subseteq\mathcal{F}_0$ *and* $\sigma\{\epsilon_t,\epsilon'_t,\epsilon''_t\}\subseteq\mathcal{F}_t\subseteq\mathcal{F}_{t+1}$, $t\ge0$.

($C_3$) *Let* $R_t=R'_t+R''_t$, $t\ge1$, *where* $R'_t=\gamma_{t+1}\nabla^Tv(\theta_t)\xi_{t+1}$, *and*

$$R''_{t+1}=\int_0^1[\nabla v(\theta_t+s(\theta_{t+1}-\theta_t))-\nabla v(\theta_t)]^T(\theta_{t+1}-\theta_t)\,ds.$$

*Then* $\sum_{t=1}^\infty\gamma_t\xi_t$ *and* $\sum_{t=1}^\infty R_t$ *converge a.s.*

PROOF.

($C_1$) Since $\mathcal{X}$ is compact, the condition ($B_2$) implies that there exists a constant $c_1\in\mathbb{R}^+$ such that

$$\|\theta_{t+1}-\theta_t\|=\|\gamma_{t+1}H(\theta_t,\mathbf{x}_{t+1})+\gamma_{t+1}^{1+\tau}\eta(\mathbf{x}_{t+1})\|\le c_1\gamma_{t+1}[V(\mathbf{x}_{t+1})+\Delta].$$

The condition ($A_1$) yields $\gamma_{t+1}/\gamma_t=O(1)$ and $|\gamma_{t+1}-\gamma_t|=O(\gamma_t\gamma_{t+1})$ for $t\to\infty$. Consequently, there exists a constant $c_2\in\mathbb{R}^+$ such that

$$\gamma_{t+1}\le c_2\gamma_t,\qquad |\gamma_{t+1}-\gamma_t|\le c_2\gamma_t^2,\qquad t\ge0.$$

Let $\epsilon_0 = \epsilon'_0 = 0$, and

$$\begin{aligned}\epsilon_{t+1} &= \gamma_{t+1}[u(\theta_t, \mathbf{x}_{t+1}) - \mathbf{P}_{\theta_t} u(\theta_t, \mathbf{x}_t)],\\ \epsilon'_{t+1} &= \gamma_{t+1}[\mathbf{P}_{\theta_{t+1}} u(\theta_{t+1}, \mathbf{x}_{t+1}) - \mathbf{P}_{\theta_t} u(\theta_t, \mathbf{x}_{t+1})]\\ &\quad + (\gamma_{t+2} - \gamma_{t+1})\mathbf{P}_{\theta_{t+1}} u(\theta_{t+1}, \mathbf{x}_{t+1}) + \gamma_{t+1}^{1+\tau}\eta(\mathbf{x}_{t+1}),\\ \epsilon''_t &= -\gamma_{t+1}\mathbf{P}_{\theta_t} u(\theta_t, \mathbf{x}_t).\end{aligned}$$

It is easy to verify that (29) is satisfied.

($C_2$) Since $\sigma(\theta_t) \subseteq \mathcal{F}_t$, we have

$$E(u(\theta_t, \mathbf{x}_{t+1})|\mathcal{F}_t) = \mathbf{P}_{\theta_t} u(\theta_t, \mathbf{x}_t),$$

which concludes (30). The condition ($B_2$) implies that there exist constants $c_3, c_4, c_5, c_6, c_7, c_8 \in \mathbb{R}^+$ and $\tau' = \min(\beta, \tau) > 0$ such that

$$\begin{aligned}\|\epsilon_{t+1}\|^2 &\le 2c_3\gamma_{t+1}^2 V^2(\mathbf{x}_t),\\ \|\epsilon'_{t+1}\| &\le c_4\gamma_{t+1} V(\mathbf{x}_{t+1})\|\theta_{t+1} - \theta_t\|^\beta + c_5\gamma_{t+1}^2 V(\mathbf{x}_{t+1}) + c_6\gamma_{t+1}^{1+\tau}\Delta\\ &\le c_7\gamma_{t+1}^{1+\tau'}[V(\mathbf{x}_{t+1}) + \Delta)],\\ \|\epsilon''_{t+1}\|^2 &\le c_8\gamma_{t+1}^2 V^2(\mathbf{x}_{t+1}).\end{aligned}$$

It follows from the condition ($A_1$) and the condition $\sup_{\mathbf{x}} V(\mathbf{x}) < \infty$ that the series $\sum_{t=0}^\infty \|\epsilon_{t+1}\|^2$, $\sum_{t=0}^\infty \|\epsilon'_t\|$ and $\sum_{t=0}^\infty \|\epsilon''_t\|^2$ all converge.

($C_3$) Let $M = \sup_{\theta\in\Theta} \max\{\|h(\theta)\|, \|\nabla v(\theta)\|\}$, and $L$ is the Lipschitz constant of $\nabla v(\cdot)$. Since $\sigma\{\theta_t\} \subset \mathcal{F}_t$, the condition ($C_2$) implies that $E(\nabla^T v(\theta_t)\epsilon_{t+1}|\mathcal{F}_t) = 0$. In addition, we have

$$\sum_{t=0}^\infty E(|\nabla^T v(\theta_t)\epsilon_{t+1}|)^2 \le M^2 \sum_{t=0}^\infty E(\|\epsilon_{t+1}\|^2) < \infty.$$

It follows from the martingale convergence theorem [Hall and Heyde (1980), Theorem 2.15] that both $\sum_{t=0}^\infty \epsilon_{t+1}$ and $\sum_{t=0}^\infty \nabla^T v(\theta_t)\epsilon_{t+1}$ converge almost surely. Since

$$\sum_{t=0}^\infty |\nabla^T v(\theta_t)\epsilon'_{t+1}| \le M \sum_{t=1}^\infty \|\epsilon'_t\|,$$

$$\sum_{t=1}^\infty \gamma_t^2\|\xi_t\|^2 \le 4\sum_{t=1}^\infty \|\epsilon_t\|^2 + 4\sum_{t=1}^\infty \|\epsilon'_t\|^2 + 8\sum_{t=0}^\infty \|\epsilon''_t\|^2,$$

it follows from ($C_2$) that both $\sum_{t=0}^\infty |\nabla^T v(\theta_t)\epsilon'_{t+1}|$ and $\sum_{t=1}^\infty \gamma_t^2\|\xi_t\|^2$ converge. In addition,

$$\begin{aligned}\|R''_{t+1}\| &\le L\|\theta_{t+1} - \theta_t\|^2 = L\|\gamma_{t+1}h(\theta_t) + \gamma_{t+1}\xi_{t+1}\|^2\\ &\le 2L(M^2\gamma_{t+1}^2 + \gamma_{t+1}^2\|\xi_{t+1}\|^2),\end{aligned}$$

$$|(\nabla v(\theta_{t+1}) - \nabla v(\theta_t))^T \epsilon''_{t+1}| \le L\|\theta_{t+1} - \theta_t\|\|\epsilon''_{t+1}\|,$$

for all $t \geq 0$. Consequently,

$$\sum_{t=1}^{\infty} |R_t''| \leq 2LM^2 \sum_{t=1}^{\infty} \gamma_t^2 + 2L \sum_{t=1}^{\infty} \gamma_t^2 \|\xi_t\|^2 < \infty,$$

$$\sum_{t=0}^{\infty} |(v(\theta_{t+1}) - v(\theta_t))^T \epsilon_{t+1}''| \leq \left(2L^2M^2 \sum_{t=1}^{\infty} \gamma_t^2 + 2L^2 \sum_{t=1}^{\infty} \gamma_t^2 \|\xi_t\|^2\right)^{1/2} \times \left(\sum_{t=1}^{\infty} \|\epsilon_t''\|^2\right)^{1/2} < \infty.$$

Since

$$\sum_{t=1}^{n} \gamma_t \xi_t = \sum_{t=1}^{n} \epsilon_t + \sum_{t=1}^{n} \epsilon_t' + \epsilon_n'' - \epsilon_0'',$$

$$\begin{aligned}\sum_{t=0}^{n} R_{t+1}' = {} & \sum_{t=0}^{n} \nabla^T v(\theta_t) \epsilon_{t+1} + \sum_{t=0}^{n} \nabla^T v(\theta_t) \epsilon_{t+1}' \\ & - \sum_{t=0}^{n} (\nabla v(\theta_{t+1} - \nabla v(\theta_t))^T \epsilon_{t+1}'' \\ & + \nabla^T v(\theta_{n+1}) \epsilon_{n+1}'' - \nabla^T v(\theta_0) \epsilon_0'',\end{aligned}$$

it is obvious that $\sum_{t=1}^{\infty} \gamma_t \xi_t$ and $\sum_{t=1}^{\infty} R_t$ converge almost surely.

The proof for Lemma A.2 is completed. □

Based on the above lemmas, Theorem A.2 can be proved in a similar way to Theorem 2.2 of Tadić (1997). Since the manuscript Tadić (1997) is not available publicly, we rewrite the proof to make the paper be self-contained.

PROOF OF THEOREM A.2. Let $M = \sup_{\theta \in \Theta} \max\{\|h(\theta)\|, |v(\theta)|\}$ and $\mathcal{V}_\varepsilon = \{\theta : v(\theta) \leq \varepsilon\}$. Applying Taylor's expansion formula [Folland (1990)], we have

$$v(\theta_{t+1}) = v(\theta_t) + \gamma_{n+1} \dot{v}(\theta_{t+1}) + R_{t+1}, \qquad t \geq 0,$$

which implies that

$$\sum_{i=0}^{t} \gamma_{i+1} \dot{v}(\theta_i) = v(\theta_{t+1}) - v(\theta_0) - \sum_{i=0}^{t} R_{i+1} \geq -2M - \sum_{i=0}^{t} R_{i+1}.$$

Since $\sum_{i=0}^{t} R_{i+1}$ converges (owing to Lemma A.2), $\sum_{i=0}^{t} \gamma_{i+1} \dot{v}(\theta_i)$ also converges. Furthermore,

$$v(\theta_t) = v(\theta_0) + \sum_{i=0}^{t-1} \gamma_{i+1} \dot{v}(\theta_i) + \sum_{i=0}^{t-1} R_{i+1}, \qquad t \geq 0,$$

$\{v(\theta_t)\}_{t\geq 0}$ also converges. On the other hand, the conditions ($A_1$) and ($A_2$) imply $\underline{\lim}_{t\to\infty} d(\theta_t,\mathcal{L})=0$. Otherwise, there exists $\varepsilon>0$ and $n_0$ such that $d(\theta_t,\mathcal{L})\geq\varepsilon$, $t\geq n_0$; as $\sum_{t=1}^{\infty}\gamma_t=\infty$ and $p=\sup\{\dot v(\theta):\theta\in\mathcal{V}_\varepsilon^c\}<0$, it is obtained that $\sum_{t=n_0}^{\infty}\gamma_{t+1}\dot v(\theta_t)\leq p\sum_{t=1}^{\infty}\gamma_{t+1}=-\infty$.

Suppose that $\overline{\lim}_{t\to\infty} d(\theta_t,\mathcal{L})>0$. Then, there exists $\varepsilon>0$ such that $\overline{\lim}_{t\to\infty} d(\theta_t,\mathcal{L})\geq 2\varepsilon$. Let $t_0=\inf\{t\geq 0: d(\theta_t,\mathcal{L})\geq 2\varepsilon\}$, while $t_k'=\inf\{t\geq t_k: d(\theta_t,\mathcal{L})\leq\varepsilon\}$ and $t_{k+1}=\inf\{t\geq t_k': d(\theta_t,\mathcal{L})\geq 2\varepsilon\}$, $k\geq 0$. Obviously, $t_k<t_{k'}<t_{k+1}$, $k\geq 0$, and

$$d(\theta_{t_k},\mathcal{L})\geq 2\varepsilon,\qquad d(\theta_{t_k'},\mathcal{L})\leq\varepsilon,\qquad d(\theta_t,\mathcal{L})\geq\varepsilon,\qquad t_k\leq t<t_k', k\geq 0.$$

Let $q=\sup\{\dot v(\theta):\theta\in\mathcal{V}_\varepsilon^c\}$. Then

$$q\sum_{k=0}^{\infty}\sum_{i=t_k}^{t_k'-1}\gamma_{i+1}\geq\sum_{k=0}^{\infty}\sum_{i=t_k}^{t_k'-1}\gamma_{i+1}\dot v(\theta_i)\geq\sum_{t=0}^{\infty}\gamma_{t+1}\dot v(\theta_t)>-\infty.$$

Therefore, $\sum_{k=0}^{\infty}\sum_{i=t_k}^{t_k'-1}\gamma_{i+1}<\infty$, and consequently, $\lim_{k\to\infty}\sum_{i=t_k}^{t_k'-1}\gamma_{i+1}=0$. Since $\sum_{t=1}^{\infty}\gamma_t\xi_t$ converges (owing to Lemma A.2), we have

$$\varepsilon\leq\|\theta_{t_k'}-\theta_{t_k}\|\leq M\sum_{i=t_k}^{t_k'-1}\gamma_{i+1}+\left\|\sum_{i=t_k}^{t_k'-1}\gamma_{i+1}\xi_{i+1}\right\|\longrightarrow 0,$$

as $k\to\infty$. This contradicts our assumption $\varepsilon>0$. Hence, $\overline{\lim}_{t\to\infty} d(\theta_t,\mathcal{L})>0$ does not hold. Therefore, $\lim_{t\to\infty} d(\theta_t,\mathcal{L})=0$ almost surely. □

**A.4. Proof of Theorem A.1.** Let $\mathbf{e}_{\mathbf{x}_t}=(e_{t1},\dots,e_{tm})$. Since the kernel used in (7) has a bounded support, $\widehat{p}_{ti}-e_{ti}/\kappa$ can then be re-expressed as

$$\widehat{p}_{ti}-e_{ti}/\kappa=\frac{\sum_{l=\max\{1,i-k_0\}}^{\min\{m,i+k_0\}}W(\Lambda l/(mh_t))(e_{t,i+l}/\kappa-e_{ti}/\kappa)}{\sum_{l=\max\{1,i-k_0\}}^{\min\{m,i+k_0\}}W(\Lambda l/(mh_t))},\tag{31}$$

where $k_0=[\frac{Cmh_t}{\Lambda}]$, and $[z]$ denotes the maximum integer less than $z$. By noting that $-1\leq\frac{e_{tj}}{\kappa}-\frac{e_{ti}}{\kappa}\leq+1$, we have $|\widehat{p}_{ti}-e_{ti}/\kappa|\leq 2k_0$. This is true even when $k_0=0$. Thus, there exists a bounded function $-2Cm/\Lambda\leq\eta_i^*(\mathbf{e}_{\mathbf{x}_t})\leq 2Cm/\Lambda$ such that

$$\widehat{p}_{ti}-e_{ti}/\kappa=h_t\eta_i^*(\mathbf{e}_{\mathbf{x}_t}).\tag{32}$$

Since $h_t$ is chosen in (9) as a power function of $\gamma_t$, the SSAMC algorithm falls into the class of stochastic approximation MCMC algorithms described in Section A.2 by letting $\eta(\mathbf{x}_t)=(\eta_1^*(\mathbf{e}_{\mathbf{x}_t}),\dots,\eta_m^*(\mathbf{e}_{\mathbf{x}_t}))$, and its convergence can be proved by verifying that it satisfies the conditions ($A_1$) to ($A_3$):

($A_1$) It is obvious that this condition is satisfied by the sequence as specified in (6).

(A$_2$) Let $\mathbf{x}_{t+1} = (x_{t+1}^{(1)}, \ldots, x_{t+1}^{(\kappa)})$, which can be regarded as a sample produced by a Markov chain on the product space $\mathcal{X}^\kappa = \mathcal{X} \times \cdots \times \mathcal{X}$ with the kernel

$$\mathbf{P}_{\theta_t}(\mathbf{x}, \mathbf{y}) = P_{\theta_t}(x^{(\kappa)}, y^{(1)}) P_{\theta_t}(y^{(1)}, y^{(2)}) \cdots P_{\theta^t}(y^{(\kappa-1)}, y^{(\kappa)}),$$

where $P_{\theta_t}(x, y)$ denotes the one-step MH kernel. To simplify notations, in the following we will drop the subscript $t$, denoting $\mathbf{x}_t$ by $\mathbf{x}$ and $\theta_t = (\theta_{t1}, \ldots, \theta_{tm})$ by $\theta = (\theta_1, \ldots, \theta_m)$.

Roberts and Tweedie (1996, Theorem 2.2) showed that if the target distribution is bounded away from 0 and $\infty$ on every compact set of its support $\mathcal{X}$, then the MH chain with a proposal distribution satisfying the local positive condition is irreducible and aperiodic, and every nonempty compact set is small. It follows from this result that $P_\theta(x, y)$ is irreducible and aperiodic, and thus $\mathbf{P}_\theta(\mathbf{x}, \mathbf{y}) = P_\theta^\kappa(x, y)$ is also irreducible and aperiodic.

Since $\mathcal{X}$ is compact, Roberts and Tweedie's result implies that $\mathcal{X}$ is a small set and the minorization condition holds on $\mathcal{X}$ for the kernel $P_\theta(x, y)$; that is, there exist an integer $l'$, a constant $\delta$ and a probability measure $\nu'(\cdot)$ such that

$$P_\theta^{l'}(x, A) \geq \delta \nu'(A) \qquad \forall x \in \mathcal{X},\ \forall A \in \mathcal{B}_{\mathcal{X}}.$$

It then follows from Rosenthal (1995, Lemma 7) that

$$\mathbf{P}_\theta^l(\mathbf{x}, \mathbf{A}) \geq \delta \nu(\mathbf{A}) \qquad \forall \mathbf{x} \in \mathcal{X}^\kappa,\ \forall \mathbf{A} \in \mathcal{B}_{\mathcal{X}^\kappa},$$

by setting $l = \min\{n : n \times \kappa \geq l', n = 1, 2, 3, \ldots\}$ and defining the measure $\nu(\cdot)$ as follows: Marginally on the first coordinate, $\nu(\cdot)$ agrees with $\nu'(\cdot)$; conditionally on the first coordinate, $\nu(\cdot)$ is defined by

$$\nu(x^{(2)}, \ldots, x^{(\kappa)} | x^{(1)}) = \mathcal{W}(x^{(2)}, \ldots, x^{(\kappa)} | x^{(1)}), \tag{33}$$

where $\mathcal{W}(x^{(2)}, \ldots, x^{(\kappa)} | x^{(1)})$ is the conditional distribution of the Markov chain samples generated by the kernel $\mathbf{P}_\theta$. Conditional on $x_t^{(1)}$, the samples $x_t^{(2)}, \ldots, x_t^{(\kappa)}$ are generated independent of all previous samples $\mathbf{x}_{t-1}, \ldots, \mathbf{x}_1$. Hence, $\mathcal{W}(x^{(2)}, \ldots, x^{(\kappa)} | x^{(1)})$ exists. This verifies condition (24) by setting $\mathbf{C} = \mathcal{X}^\kappa$. Thus, for any $\theta \in \Theta$ the following conditions hold:

$$\begin{aligned} &\mathbf{P}_\theta^l V^\alpha(\mathbf{x}) \leq \lambda V^\alpha(\mathbf{x}) + b I(\mathbf{x} \in \mathbf{C}) \qquad \forall \mathbf{x} \in \mathcal{X}^\kappa, \\ &\mathbf{P}_\theta V^\alpha(\mathbf{x}) \leq \varsigma V^\alpha(\mathbf{x}) \qquad \forall \mathbf{x} \in \mathcal{X}^\kappa, \end{aligned} \tag{34}$$

by choosing $V(\mathbf{x}) = 1$, $0 < \lambda < 1$, $b = 1 - \lambda$, $\varsigma > 1$, and $\alpha \geq 2$. These conclude that (A$_2$)(i) is satisfied.

Let $H^{(i)}(\theta, \mathbf{x})$ be the $i$th component of the vector $H(\theta, \mathbf{x}) = (\mathbf{e}_\mathbf{x}/\kappa - \boldsymbol{\pi})$. By construction, $|H^{(i)}(\theta, \mathbf{x})| = |e_\mathbf{x}^{(i)}/\kappa - \pi_i| < 1$ for all $\mathbf{x} \in \mathcal{X}^\kappa$ and $i = 1, \ldots, m$.

Therefore, there exists a constant $c_1 = \sqrt{m}$ such that for any $\theta \in \Theta$ and all $\mathbf{x} \in \mathcal{X}^\kappa$,

$$\|H(\theta, \mathbf{x})\| \le c_1. \tag{35}$$

Also, $H(\theta, \mathbf{x})$ does not depend on $\theta$ for a given sample $\mathbf{x}$. Hence, $H(\theta, \mathbf{x}) - H(\theta', \mathbf{x}) = 0$ for all $(\theta, \theta') \in \Theta \times \Theta$, and the following condition holds for the SSAMC algorithm:

$$\|H(\theta, \mathbf{x}) - H(\theta', \mathbf{x})\| \le c_1 \|\theta - \theta'\|, \tag{36}$$

for all $(\theta, \theta') \in \Theta \times \Theta$. Equations (35) and (36) imply that $(\mathrm{A}_2)$(ii) is satisfied by choosing $\beta = 1$ and $V(\mathbf{x}) = 1$.

Let $s_\theta(x, y) = q(x, y) \min\{1, r_\theta(x, y)\}$, where $r_\theta(x, y) = \frac{f_\theta(y) q(y,x)}{f_\theta(x) q(x,y)}$. Thus, we have

$$\begin{aligned}\left|\frac{\partial s_\theta(x, y)}{\partial \theta_i}\right| &= |-q(x, y) I(r_\theta(x, y) < 1) \\ &\qquad \times I(J(x) = i \text{ or } J(y) = i) I(J(x) \ne J(y)) r_\theta(x, y)| \\ &\le q(x, y),\end{aligned}$$

where $I(\cdot)$ is the indicator function, and $J(x)$ denotes the index of the subregion to which $x$ belongs. The mean-value theorem implies that there exists a constant $c_2$ such that

$$|s_\theta(x, y) - s_{\theta'}(x, y)| \le q(x, y) c_2 \|\theta - \theta'\|, \tag{37}$$

which implies that

$$\sup_x \int_{\mathcal{X}} |s_\theta(x, y) - s_{\theta'}(x, y)| \, dy \le c_2 \|\theta - \theta'\|. \tag{38}$$

Since the MH kernel can be expressed in the form

$$P_\theta(x, dy) = s_\theta(x, dy) + I(x \in dy) \left[1 - \int_{\mathcal{X}} s_\theta(x, z) \, dz\right],$$

for any measurable set $A \subset \mathcal{X}$ we have

$$\begin{aligned}&|P_\theta(x, A) - P_{\theta'}(x, A)| \\ &\qquad = \left| \int_A \left[s_\theta(x, y) - s_{\theta'}(x, y)\right] dy \right. \\ &\qquad\qquad \left. + I(x \in A) \int_{\mathcal{X}} \left[s_{\theta'}(x, z) - s_\theta(x, z)\right] dz \right| \\ &\qquad \le \int_{\mathcal{X}} |s_\theta(x, y) - s_{\theta'}(x, y)| \, dy + I(x \in A) \int_{\mathcal{X}} |s_{\theta'}(x, z) - s_\theta(x, z)| \, dz \\ &\qquad \le 2 \int_{\mathcal{X}} |s_\theta(x, y) - s_{\theta'}(x, y)| \, dy \\ &\qquad \le 2 c_2 \|\theta - \theta'\|.\end{aligned} \tag{39}$$

Since $\mathbf{P}_\theta(\mathbf{x},\mathbf{A})$ can be expressed in the following form:

$$\begin{aligned}\mathbf{P}_\theta(\mathbf{x},\mathbf{A})\\ &=\int_{A_1}\cdots\int_{A_\kappa} P_\theta(x^{(\kappa)},y^{(1)})P_\theta(y^{(1)},y^{(2)})\cdots P_\theta(y^{(\kappa-1)},y^{(\kappa)})\,dy^{(1)}\cdots dy^{(\kappa)},\end{aligned}$$

(39) implies that there exists a constant $c_3$ such that

$$\begin{aligned}&|\mathbf{P}_\theta(\mathbf{x},\mathbf{A})-\mathbf{P}_{\theta'}(\mathbf{x},\mathbf{A})|\\ &\quad=\bigg|\int_{A_1}\cdots\int_{A_\kappa}[P_\theta(x^{(\kappa)},y^{(1)})\\ &\qquad\qquad\times P_\theta(y^{(1)},y^{(2)})\cdots P_\theta(y^{(\kappa-1)},y^{(\kappa)})\\ &\qquad\qquad-P_{\theta'}(x^{(\kappa)},y^{(1)})P_{\theta'}(y^{(1)},y^{(2)})\cdots\\ &\qquad\qquad\qquad\times P_{\theta'}(y^{(\kappa-1)},y^{(\kappa)})]\,dy^{(1)}\cdots dy^{(\kappa)}\bigg|\\ &\quad\le\int_{A_1}\int_{\mathcal{X}}\cdots\int_{\mathcal{X}}|P_\theta(x^{(\kappa)},y^{(1)})-P_{\theta'}(x^{(\kappa)},y^{(1)})|\\ &\qquad\qquad\times P_\theta(y^{(1)},y^{(2)})\cdots P_\theta(y^{(\kappa-1)},y^{(\kappa)})\,dy^{(1)}\cdots dy^{(\kappa)}\\ &\qquad+\int_{\mathcal{X}}\int_{A_2}\int_{\mathcal{X}}\cdots\int_{\mathcal{X}}P_{\theta'}(x^{(\kappa)},y^{(1)})|P_\theta(y^{(1)},y^{(2)})-P_{\theta'}(y^{(1)},y^{(2)})|\\ &\qquad\qquad\times P_\theta(y^{(2)},y^{(3)})\cdots P_\theta(y^{(\kappa-1)},y^{(\kappa)})\,dy^{(1)}\cdots dy^{(\kappa)}\\ &\qquad+\cdots\\ &\qquad+\int_{\mathcal{X}}\cdots\int_{\mathcal{X}}\int_{A_\kappa}P_{\theta'}(x^{(\kappa)},y^{(1)})\cdots P_{\theta'}(y^{(\kappa-2)},y^{(\kappa-1)})\\ &\qquad\qquad\times|P_\theta(y^{(\kappa-1)},y^{(\kappa)})-P_{\theta'}(y^{(\kappa-1)},y^{(\kappa)})|\,dy^{(1)}\cdots dy^{(\kappa)}\\ &\quad\le c_3\|\theta-\theta'\|,\end{aligned}$$

which implies that (27) is satisfied.

For any function $g\in\mathcal{L}_V$,

$$\begin{aligned}\|\mathbf{P}_\theta g-\mathbf{P}_{\theta'}g\|_V&=\left\|\int(\mathbf{P}_\theta(\mathbf{x},d\mathbf{y})-\mathbf{P}_{\theta'}(\mathbf{x},d\mathbf{y}))g(\mathbf{y})\right\|_V\\ &=\bigg\|\int_{\mathcal{X}_+^\kappa}(\mathbf{P}_\theta(\mathbf{x},d\mathbf{y})-\mathbf{P}_{\theta'}(\mathbf{x},d\mathbf{y}))g(\mathbf{y})\\ &\qquad+\int_{\mathcal{X}_-^\kappa}(\mathbf{P}_\theta(\mathbf{x},d\mathbf{y})-\mathbf{P}_{\theta'}(\mathbf{x},d\mathbf{y}))g(\mathbf{y})\bigg\|_V\\ &\le\|g\|_V\{|\mathbf{P}_\theta(\mathbf{x},\mathcal{X}_+^\kappa)-\mathbf{P}_{\theta'}(\mathbf{x},\mathcal{X}_+^\kappa)|+|\mathbf{P}_\theta(\mathbf{x},\mathcal{X}_-^\kappa)-\mathbf{P}_{\theta'}(\mathbf{x},\mathcal{X}_-^\kappa)|\}\\ &\le 4c_2\|g\|_V|\theta-\theta'|\qquad[\text{following from (39)}]\end{aligned}$$

where $\mathcal{X}_+^\kappa = \{\mathbf{y} : \mathbf{y} \in \mathcal{X}^\kappa, \mathbf{P}_\theta(\mathbf{x}, d\mathbf{y}) - \mathbf{P}_{\theta'}(\mathbf{x}, \mathbf{y}) > 0\}$ and $\mathcal{X}_-^\kappa = \mathcal{X}^\kappa \setminus \mathcal{X}_+^\kappa$. Therefore, condition ($A_2$)(iii) is satisfied by choosing $V(\mathbf{x}) = 1$ and $\beta = 1$.

($A_3$) Since the invariant distribution of the kernel $P_\theta(x, \cdot)$ is $f_\theta(x)$, we have for any fixed $\theta$,

$$E(e_{\mathbf{x}}^{(i)}/\kappa - \pi_i) = \frac{\int_{E_i} \psi(x)\,dx/e^{\theta_i}}{\sum_{k=1}^m [\int_{E_k} \psi(x)\,dx/e^{\theta_k}]} - \pi_i = \frac{S_i}{S} - \pi_i, \qquad i = 1, \dots, m, \tag{40}$$

where $S_i = \int_{E_i} \psi(x)\,dx/e^{\theta_i}$ and $S = \sum_{k=1}^m S_k$. Thus, we have

$$h(\theta) = \int_{\mathcal{X}} H(\theta, \mathbf{x}) f(d\mathbf{x}) = \left(\frac{S_1}{S} - \pi_1, \dots, \frac{S_m}{S} - \pi_m\right)'.$$

It follows from (40) that $h(\theta)$ is a continuous function of $\theta$. Let $v(\theta) = \frac{1}{2}\sum_{k=1}^m (\frac{S_k}{S} - \pi_k)^2$. As shown below, $v(\theta)$ has continuous partial derivatives of the first order.

Solving the system of equations formed by (40), we have

$$\mathcal{L} = \left\{(\theta_1, \dots, \theta_m) : \theta_i = \text{Const} + \log\left(\int_{E_i} \psi(\mathbf{x})\,d\mathbf{x}\right) - \log(\pi_i), i = 1, \dots, m; \theta \in \Theta\right\},$$

where $\text{Const} = \log(S)$ can be determined by imposing a constraint on $S$. For example, setting $S = 1$ leads to that $c = 0$. It is obvious that $\mathcal{L}$ is nonempty and $v(\theta) = 0$ for every $\theta \in \mathcal{L}$.

To verify the conditions related to $\dot{v}(\theta)$, we have the following calculations:

$$\begin{aligned} \frac{\partial S}{\partial \theta_i} &= \frac{\partial S_i}{\partial \theta_i} = -S_i, \\ \frac{\partial S_i}{\partial \theta_j} &= \frac{\partial S_j}{\partial \theta_i} = 0, \\ \frac{\partial (S_i/S)}{\partial \theta_i} &= -\frac{S_i}{S}\left(1 - \frac{S_i}{S}\right), \\ \frac{\partial (S_i/S)}{\partial \theta_j} &= \frac{\partial (S_j/S)}{\partial \theta_i} = \frac{S_i S_j}{S^2}, \end{aligned} \tag{41}$$

for $i, j = 1, \ldots, m$ and $i \neq j$,

$$
\begin{aligned}
\frac{\partial v(\theta)}{\partial \theta_i} &= \frac{1}{2} \sum_{k=1}^{m} \frac{\partial (S_k/S - \pi_k)^2}{\partial \theta_i} \\
&= \sum_{j \neq i} \left( \frac{S_j}{S} - \pi_j \right) \frac{S_i S_j}{S^2} - \left( \frac{S_i}{S} - \pi_i \right) \frac{S_i}{S} \left( 1 - \frac{S_i}{S} \right) \\
&= \sum_{j=1}^{m} \left( \frac{S_j}{S} - \pi_j \right) \frac{S_i S_j}{S^2} - \left( \frac{S_i}{S} - \pi_i \right) \frac{S_i}{S} \\
&= \mu_{\eta^*} \frac{S_i}{S} - \left( \frac{S_i}{S} - \pi_i \right) \frac{S_i}{S},
\end{aligned} \tag{42}
$$

for $i = 1, \ldots, m$, where $\mu_{\eta^*} = \sum_{j=1}^{m} (\frac{S_j}{S} - \pi_j) \frac{S_j}{S}$. Thus, we have

$$
\begin{aligned}
\dot{v}(\theta) &= \mu_{\eta^*} \sum_{i=1}^{m} \left( \frac{S_i}{S} - \pi_i \right) \frac{S_i}{S} - \sum_{i=1}^{m} \left( \frac{S_i}{S} - \pi_i \right)^2 \frac{S_i}{S} \\
&= -\left\{ \sum_{i=1}^{m} \left( \frac{S_i}{S} - \pi_i \right)^2 \frac{S_i}{S} - \mu_{\eta^*}^2 \right\} \\
&= -\sigma_{\eta^*}^2 \leq 0,
\end{aligned} \tag{43}
$$

where $\sigma_{\eta^*}^2$ denotes the variance of the discrete distribution defined in the following table:

| | | | |
|---|---|---|---|
| State ($\eta^*$) | $\frac{S_1}{S} - \pi_1$ | $\cdots$ | $\frac{S_m}{S} - \pi_m$ |
| Prob. | $\frac{S_1}{S}$ | $\cdots$ | $\frac{S_m}{S}$ |

If $\theta \in \mathcal{L}$, $\dot{v}(\theta) = 0$; otherwise, $\dot{v}(\theta) < 0$. Therefore, $\sup_{\theta \in Q} \dot{v}(\theta) < 0$ for any compact set $Q \subset \mathcal{L}^c$.

The proof is completed.

**Acknowledgments.** The author thanks Dr. Vladislav Tadić for sending him his manuscript on the convergence of stochastic iterative algorithms, thanks Dr. Chuanhai Liu for his early discussion of this project, and thanks the Editor, Associate Editor and two referees for their constructive comments which have led to significant improvement of this paper.

Department of Statistics
Texas A&M University
College Station, Texas 77843-3143
USA
E-mail: fliang@stat.tamu.edu